\documentclass[journal]{IEEEtran}

\pdfoutput=1

\usepackage{cite}

\ifCLASSINFOpdf
   \usepackage[pdftex]{graphicx}
\else
\fi
\usepackage{amssymb}
\usepackage{amsmath}
\usepackage{dcolumn}
\newcolumntype{d}[2]{D{.}{.}{#1.#2}}
\newcolumntype{.}{D{.}{.}{-1}}
\usepackage{multirow}
\usepackage{colortbl}
\usepackage{array}

\ifCLASSOPTIONcompsoc
  \usepackage[caption=false,font=normalsize,labelfont=sf,textfont=sf]{subfig}
\else
  \usepackage[caption=false,font=footnotesize]{subfig}
\fi

\usepackage{fixltx2e}
\usepackage{url}

\newcommand{\ie}{\emph{i.e.}}
\newcommand{\eg}{\emph{e.g.}}

\newcommand{\argmin}{\mathop{\mathrm{argmin}}}
\newcommand{\prob}{\mathrm{P}}
\newcommand{\diag}{\mathrm{diag}}
\newcommand{\blkdiag}{\mathrm{blkdiag}}
\newcommand{\E}{\mathsf{E}}
\newcommand{\ones}{\mathbf{1}}
\newcommand{\rv}[1]{\mathbf{#1}}
\newcommand{\reals}{\mathbb{R}}
\renewcommand{\vec}{\mathrm{vec}}
\newcommand{\prox}{\mathrm{prox}}

\newcommand{\vml}{\hat v_{\mathrm{f}}}
\newcommand{\vmly}{\hat v_{\mathrm{y}}}
\newcommand{\uml}{\hat u_{\mathrm{y}}}

\usepackage{tikz,pgfplots}
\pgfplotsset{compat=1.13}
\usepgfplotslibrary{groupplots}
\usetikzlibrary{matrix,spy}

\newcommand{\remove}[1]{}
\newcommand{\add}[1]{#1}

\begin{document}

\title{A Convex Reconstruction Model for X-ray Tomographic Imaging with Uncertain Flat-fields}

\author{Hari~Om~Aggrawal,~Martin~S.~Andersen, Sean Rose, and Emil Y. Sidky%
\thanks{Hari Om Aggrawal and Martin S. Andersen are with the
  Department of Applied Mathematics and Computer Science, Technical
  University of Denmark,  2800 Lyngby, Denmark (e-mail: \{haom,
  mskan\}@dtu.dk).}
\thanks{Sean Rose and Emil Y. Sidky are with the Department of
  Radiology, University of Chicago, 5841 South Maryland Avenue,
  Chicago, Illinois 60637, USA (email: seanrose949@gmail.com, sidky@uchicago.edu).}%
}

\maketitle

\begin{abstract}
  Classical methods for X-ray computed tomography are based on the
  assumption that the X-ray source intensity is known, but in
  practice, the intensity is measured and hence uncertain. Under
  normal operating conditions, when the exposure time is sufficiently
  high, this kind of uncertainty typically has a negligible effect on
  the reconstruction quality. However, in time- or dose-limited
  applications such as dynamic CT, this uncertainty may cause severe
  and systematic artifacts known as ring artifacts. By carefully
  modeling the measurement process and by taking uncertainties into
  account, we derive a new convex model that leads to improved
  reconstructions despite poor quality measurements.  We demonstrate
  the effectiveness of the methodology based on simulated and real
  data sets.
\end{abstract}


%
\IEEEpeerreviewmaketitle

\section{Introduction}
\IEEEPARstart{X}{-ray} computed tomography (CT) is a non-invasive
method that is used to image the internal structure of objects without
cutting or breaking them. An X-ray source illuminates an object from
different directions while detectors capture the attenuated X-rays. As
the X-rays propagate through the object along straight lines, they are
attenuated exponentially with a rate of decay that depends on the
material. This relationship is explained by the Lambert--Beer law which
forms the basis of major X-ray CT reconstruction models and methods;
see \eg\ \cite{Buzug}. Reconstruction methods
estimate the spacial attenuation of the object of interest based on a
number of X-ray images, given the measurement geometry, the
source intensity, and possibly some assumptions on the statistical
nature of the measurement process.

In practice, the source intensity is never known exactly, but it is
estimated by acquiring a number of X-ray images without an object in
the scanner. Such measurements are also known as air scans
\cite{Whiting2006}, flat-fields, or white-fields
\cite{Sijbers2004}. The elementwise mean of these measurements
provides an estimate of the flat-field intensity and may be used for
computing reconstructions. However, in practice the measurements are
noisy, and hence the flat-field intensity estimate is a random
variable whose variance is proportional to the ratio of the flat-field
intensity and number of flat-field samples
\cite{MorrisH.DeGroot}. Consequently, the signal-to-noise ratio (SNR)
of the flat-field intensity estimate is proportional to the square
root of the product of the flat-field intensity and the number of
samples. Therefore, if the flat-field intensity is low or if the
number of flat-field measurements is small, the flat-field estimation
error may be significant and lead to reconstruction artifacts and
errors. Since the flat-field estimate is used to normalize
measurements from all projection directions, the estimation errors
result in systematic reconstruction errors. These are known as ring
artifacts\cite{Kowalski1977} since they appear as concentric circles
superimposed on the reconstruction, and they are a common problem that
can mask important features in the reconstructed image
\cite{Thomas2010,Dahlman2012}. Ring artifacts may not only occur
because of flat-field estimation errors; miscalibrated or dead
detector elements and non-uniform sensitivities may also systematically
corrupt the measurements and lead to ring artifacts in the
reconstruction \cite{Sijbers2004}.

An experimental study \cite{Fahrig2006} has pointed out that the ring
artifacts are more severe when the X-ray source intensity is low, and
hence a reconstruction from low-intensity measurements may be very
sensitive to the assumptions upon which the reconstruction method is
based. The problem may arise when the acquisition time is limited,
\eg, in dynamic or time-resolved tomography, or if the application
imposes strict dose limitations. Thus, tomographic reconstruction based on
low-intensity measurements is a challenging problem, in part because of the low SNR.

One approach to combating ring artifacts is to move the detector array
between projections \cite{DKB+:01}. This has an averaging effect on
the systematic error due to flat-field estimation errors and often
results in noticeable improvements, but it does not address or model
the underlying cause. Moreover, it requires special hardware for the
acquisition, and it is not suited for applications such as dynamic CT
where fast acquisition times are important. Alternative software-based
methods to mitigate ring artifacts also exist. Roughly speaking, these
methods can be put into three categories: sinogram preprocessing
methods \cite{Kowalski1978,Raven1998,Munch2009,Rashid2012,Kim2014},
combined ring reduction and reconstruction methods
\cite{Paleo2015,AdityaMohan2015}, and post-processing methods that
reduce or remove rings from a reconstruction
\cite{Sijbers2004,Prell2009,Yan2016}. The preprocessing
methods detect and remove/reduce stripes in the sinogram which, in turn, reduces
the ring artifacts in the image domain. These algorithms are typically
based on Fourier domain filtering \cite{Raven1998}, wavelet domain
filtering\cite{Munch2009}, or a normalization of measurements by
estimating the sensitivity of each detector pixel \cite{Kim2014}.  The
post-processing methods transform the reconstructed image from
Cartesian to polar coordinates \cite{Sijbers2004} and remove stripes
using, \eg, a median filter\cite{Prell2009}, a wavelet filter, or a
variational model for destriping \cite{Yan2016}.

In two recently proposed methods \cite{AdityaMohan2015,Paleo2015},
ring artifact correction is included as an intrinsic part of the
reconstruction process. Motivated by the cause of ring artifacts,
which appear as stripes in the sinogram domain, the sinogram is split
into the sum of the true sinogram and a component which represents the
systematic stripe errors.  Although the combined ring-reduction and
reconstruction methods do take the systematic nature of the flat-field
estimation errors in the sinogram domain into account, they do not
explicitly model the source of the errors nor their statistical
properties.

Existing methods for mitigating ring artifacts have been shown to work
reasonably well when applied to measurement data with high or
acceptable SNRs. However, we are not aware of any studies that
investigate ring artifact correction for low SNR measurements and
where the intensity of X-ray beam is assumed to be uncertain. To this
end, we derive a new reconstruction model that is based on a rigorous
statistical description of our model assumptions.  Unlike existing
correction methods that, roughly speaking, are based on the geometric
nature of ring artifacts in either the sinogram or the reconstruction,
our approach is based on a model of a fundamental cause of these
artifacts. The resulting reconstruction method jointly estimates the
flat-field and the attenuation image, and we show that the estimation
problem can be solved efficiently by solving a convex optimization
problem. We also derive a quadratic approximation model which is
similar to an existing weighted least-squares reconstruction model.

\subsubsection*{Outline} Section \ref{conventional} introduces our
model assumptions and reviews some existing approaches to CT
reconstruction based on low SNR measurements. We illustrate the
sensitivity of these existing methods to flat-field intensity
estimation errors. Section \ref{proposed} proposes a new
reconstruction model and discusses different parameter selection
strategies. We describe our numerical implementation in Section
\ref{s-implementation}, and we validate the proposed model based on
simulated data as well as real tomographic measurements in Section
\ref{results}.  Section \ref{conclusion} concludes the paper.

\subsubsection*{Notation} The set $\reals^n$ denotes the
$n$-dimensional real space, $\reals^n_+$ is the nonnegative orthant of
$\reals^n$, and $\reals^{m\times n}$ is the set of $m\times n$
real-valued matrices. Upper case letters denote matrices, lower case
letters denote vectors or scalars, and boldface letters denote random
variables. Given a vector $x \in \reals^n$, the matrix $\diag(x)$ is
the $n\times n$ diagonal matrix with the elements of $x$ on the
diagonal. Similarly, given a set of $r$ square matrices
$S_1,\ldots,S_r$, the matrix $\blkdiag(S_1,\ldots,S_r)$ denotes the
block-diagonal matrix with diagonal blocks $S_1,\ldots,S_r$. The
vector $e_i$ denotes the $i$th column of an identity matrix, and
$\ones$ denotes a vector of ones. Given a vector $x\in \reals^n$, the
notation $\log(x)$ and $\exp(x)$ is interpreted as elementwise
logarithm and exponentiation\add{. $A \otimes B$ denotes
  the Kronecker product of $A \in \reals^{m\times n}$ and $B \in
  \reals^{p \times q}$, $\|A\|_F$ denotes the Frobenius norm of $A$, and $|A| \in \reals^{m\times n}$ is the element-wise absolute value of $A$.} The vector $y= \vec(Y)$ denotes the
vector obtained by stacking the columns of the matrix $Y$. Given a
discrete random variable $\rv{y}$, the probability of $\rv{y} = y$ is
$\prob(\rv{y} = y)$, or using shorthand notation,
$\prob(y)$. Similarly, given a continuous random variable $\rv{z}$,
$\prob(z)$ is shorthand for the probability density associated with
$\rv{z}$, evaluated at $z$, and finally, $\E[\rv{z}]$ denotes the
expectation of $\rv{z}$.

\section{Conventional Reconstruction Approach} \label{conventional}

\subsection{System and Measurement Model}\label{ss-model}

The Lambert--Beer law describes how an X-ray beam is attenuated as it
travels through an object that is characterized by a spatial
attenuation function $\mu(x)$. Specifically, the incident intensity of
an X-ray beam on a detector is given by
\begin{align}
  I \approx I_0 \exp\left( -\int_l \mu(x) \, dx \right)
\end{align}
where $I_0$ is the intensity of the X-ray source, and $l$ denotes the
line segment between the source and a detector. This description
does not take the detector efficiency and the statistical nature of the photon arrival process
into account.  For photon-counting detectors, it is common to assume
that the photon arrival process is a Poisson process, and
each measurement is assumed to be a sample from a Poisson distribution
whose mean is prescribed by the Lambert--Beer law.
Here we will consider a two-dimensional geometry
where $p$ projections are acquired using a one-dimensional detector array
with $r$ detector elements. We will use the notation $y_{ij}$ to
denote the measurement obtained with detector
element $i$ and projection $j$, and we will assume that the $i$th
detector element has efficiency $\eta_i \in (0,1]$ such that the
effective intensity is $v_i = \eta_i I_0$. Thus, with the assumption that the
arrival process is Poisson process, $y_{ij}$ is a
realization of a random variable $\rv{y}_{ij}$ which, conditioned on
$\mu$ and $v_i$, is a Poisson random variable whose mean is prescribed by the
Lambert--Beer law, \ie,
\begin{align}\label{e-meas-model}
  \rv{y}_{ij} \mid \mu, v_i \sim \mathrm{Poisson}\left( v_i \exp\left(
  -\int_{l_{ij}} \mu(x) \, dx \right) \right)
\end{align}
where $l_{ij}$ notes the line segment between the $i$th detector element and the
source for projection angle $j$.
For ease of notation, we define a matrix random variable $\rv{Y}$ of size
$r \times p$ with
elements $\rv{y}_{ij}$, and similarly, the $r \times p$ matrix $Y$
denotes a realization of $\rv{Y}$ and $y = \vec(Y)$.

The attenuation function $\mu(x)$ may be discretized by using a parameterization
\begin{align}\label{e-uparam}
	\mu(x) = \sum_{k=1}^{n} u_k \mu_k(x)
\end{align}
where $\mu_k(x)$ is one of $n$ basis functions (\emph{e.g.}, a pixel or voxel basis), and
$u \in \reals^n$ is a vector of unknowns (\emph{e.g.}, pixel or
voxel values). With this parameterization, the line integrals in
\eqref{e-meas-model} can be
expressed as
\[ \int_{l_{ij}} \mu(x) \, dx =  e_i^TA_ju \]
where the elements of the matrix $A_j \in \reals^{r \times n}$ are given by
\begin{equation*}
  (A_j)_{ik} = \int_{l_{ij}} \mu_k(x) \, dx,
\end{equation*}
and hence the columns of $\rv{Y}$ satisfy
\begin{align*}
\E[\rv{y}_j| u,v] = \diag(v) \exp(-A_ju), \ j=1,\ldots,p
\end{align*}
where $v = (v_1,\ldots,v_r)$.

In practice, the vector $v$ is unknown and must be measured. As
mentioned in the introduction, the measurements of $v$ are often
referred to as flat-field measurements and are simply measurements
obtained without any object in the CT scanner. We will assume that $s$
flat-field measurements are acquired for each detector element based on
the flat-field measuring model
\begin{align}\label{e-meas-model-ff}
  \rv{f}_{ij} \mid v_i \sim \mathrm{Poisson}\left( v_i \right)
\end{align}
for $i=1,\ldots,r$ and $j=1,\ldots,s$, and $\rv{F}$ denotes a
$r \times s$ matrix random variable with elements $\rv{f}_{ij}$. As
for the measurements $Y$, the matrix $F \in
\reals^{r \times s}$ denotes a realization of $\rv{F}$.

\subsection{Maximum Likelihood Estimation}

Given the flat-field measurements $F$, a maximum likelihood (ML) estimate of $v$
is given by
\begin{align}
\label{empmeanest}
  \vml &= \argmin_{v} \left\{ -\log \prob( F \mid v) \right\} \\
\notag
  &= \argmin_{v} \left\{ s \ones^Tv - \ones^T F^T \log(v)
    \right\} = \frac{1}{s} F \ones,
\end{align}
\ie, $\vml$ is simply the arithmetic average of the $s$ flat-field
measurements. This estimate can be used to compute an \emph{approximate} ML
estimate of the vector $u$ which is given by
\begin{align} \label{e-ml-aml}
  \uml &= \argmin_{u} \left\{ -\log
  \prob(Y \mid u, \vml )\right\} \\
\notag  &=\argmin_{u}  \left\{ (\ones \otimes \vml)^T \exp(-Au) +
    y^TAu\right\}
\end{align}
where $A \in \reals^{rp \times n}$ is defined as
$ A = [A_1^T \ \cdots \ A_p^T]^T$. The estimation problem
\eqref{e-ml-aml} is a convex optimization problem, and it is
essentially an approximate ML estimation problem since with our model
assumptions, the true likelihood $\prob(Y \mid u, v)$ is a function of
both $u$ and $v$. We will return to this issue in the next section.

If $y$ is positive, a quadratic approximation of \eqref{e-ml-aml} can be obtained by means
of a second-order Taylor expansion of the likelihood function
\cite{SaB:93}, and this yields the following weighted least-squares
objective function
\begin{align}\label{e-ml-wls-cost}
 \frac{1}{2} \|
  \diag(y)^{1/2} (Au- b)\|_2^2
\end{align}
where $ b = \ones \otimes \log(\vml) -
\log(y)$. Notice that if $A$ has full rank and $rp \leq n$, both
\eqref{e-ml-aml} and the quadratic approximation
\eqref{e-ml-wls-cost} reduce to the problem of solving the consistent system of
equations $Au = b$, but the two problems are generally different
when the system of equations $Au=b$ is inconsistent. The noise
properties of reconstructions based on the weighted least-squares
objective \eqref{e-ml-wls-cost} have been studied in \cite{Rose2015}.

\subsection{The Effect of Flat-field Estimation Errors} \label{effecffferror}

The flat-field estimate $\vml$ in (\ref{empmeanest}) satisfies
$\E[\vml] = v$, and hence it is an unbiased estimate. However, $\vml$
is itself a random variable with covariance $(1/s) \diag(v)$, and the
flat-field estimation error may lead to artifacts in the
reconstruction. To study how flat-field estimation errors influence
the reconstruction, we now consider a simplified model based on
Gaussian approximations. Specifically, we assume that
$(\vml)_i | v_i \sim \mathcal{N}(v_i, s^{-1} v_i)$ and
$\rv{y}_{ij} | v_i, u \sim \mathcal{N}(v_i \exp(-e_i^TA_ju), v_i
\exp(-e_i^TA_ju))$.
With these assumptions,
$ \rv{b}_{ij} = \log( (\vml)_i) - \log(\rv{y}_{ij})$ can be
approximated by linearizing each of the log terms around the mean of
their arguments, \ie,
\begin{align*}
\rv{b}_{ij} & \approx
\log(v_i) +\frac{(\vml)_i  - v_i}{v_i} -
  \log(\E[\rv{y}_{ij}]) - \frac{\rv{y}_{ij}  -
  \E[\rv{y}_{ij}]}{\E[\rv{y}_{ij}]} \\
&= e_i^TA_ju + \rv{z}_i + \rv{w}_{ij}
\end{align*}
for $i=1,\ldots,r$ and $j=1,\ldots,p$, and
where \[\mathbf{z}_i = ((\vml)_i  - v_i)/v_i, \ \rv{z}_i \sim
\mathcal{N}(0,(s v_i)^{-1})\]
and
\[ \rv{w}_{ij} = (\rv{y}_{ij} -
\E[\rv{y}_{ij}])/\E[\rv{y}_{ij}], \  \rv{w}_{ij} \sim \mathcal{N}(0,
v_i^{-1} \exp(e_i^TA_ju)).\] The terms $\rv{z}_i$
arise because of the flat-field estimation errors, and the terms $\rv{w}_{ij}$ represent the effect of measurement noise. If we define $\rv{z} =
(\rv{z}_1,\ldots,\rv{z}_r)$ and $\rv{w} = \vec(\rv{W})$ where $\rv{W}$
is the $r \times p$ matrix with elements $\rv{w}_{ij}$, then
\begin{align}
  \label{e-linear-approximation}
  \rv{b}\approx Au + \ones \otimes \rv{z} + \rv{w}.
\end{align}
Not surprisingly, this shows that flat-field estimation errors affect
all projections, and hence give rise to structured errors\remove{, and the
linear approximation reaffirms that the variance of the flat-field
errors is inversely proportional to the flat-field intensity}.

\add{
The linear approximation reaffirms that the variance of the flat-field
errors is inversely proportional to the flat-field intensity and the
number of flat-field measurements $s$. Thus, if $s$ is sufficiently
large, the flat-field estimation errors play a negligible
role. However, a twofold reduction of the flat-field error-to-noise ratio 
\[ \sqrt{ \frac{\E[\rv{z}_i^2]}{ \E[\rv{w}_{ij}^2]} } =
  \frac{1}{\sqrt{s \exp(-e_i^TA_ju)}} \] requires a fourfold increase
in the number of flat-field samples, and hence it may require many
samples to obtain a sufficiently small flat-field error-to-noise ratio.}

We now demonstrate the effect of flat-field estimation errors by
considering the behavior of reconstructions based on
\eqref{e-ml-aml}. We will use a constant flat-field $v = \omega \ones$
for $\omega > 0$ to generate a set of measurements according to the
model \eqref{e-meas-model} with $r=200$ detector elements and $p=720$
parallel beam projections covering a full rotation. \remove{For the
reconstruction we use an estimate $\vml$ where only three elements of
$\vml$ contain an estimation error of magnitude $2\sqrt{\omega}$,
corresponding to two times the standard deviation of the ML flat-field
estimate based on \eqref{empmeanest}.}
\add{
For the reconstruction we use the flat-field ML estimate $\vml$, as
defined in \eqref{empmeanest}, where only one flat-field sample ($s=1$) is acquired for each detector element based on \eqref{e-meas-model-ff}.
}

\begin{figure}
\centering
\subfloat[Phantom]{\includegraphics[width=0.45\columnwidth]{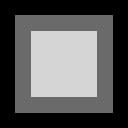}\label{acadphantom}}
\hfil
\subfloat[$\omega = 10^3$ photons]{\includegraphics[width=0.45\columnwidth]{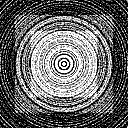}\label{acad10}}

\subfloat[$\omega = 10^4$ photons]{\includegraphics[width=0.45\columnwidth]{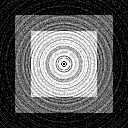}\label{acad103}}
\hfil
\subfloat[$\omega =10^5$ photons]{\includegraphics[width=0.45\columnwidth]{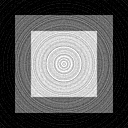}\label{acad105}}
\caption{Phantom (a) and reconstructions (b), (c), and (d), based on
  \eqref{e-ml-aml} with \remove{three} flat-field estimation errors. The
  display range for each of the images is $[0,0.6]$.}
\label{acadex}
\end{figure}

Our object $u$, shown in Fig.\ \ref{acadphantom}, consists of three
squares of different sizes where the attenuation of the innermost
square is 0.5 cm$^{-1}$, the enclosing square has attenuation
0.25 cm$^{-1}$, and the outermost square has no attenuation. The
domain size is 1 cm, and the reconstruction grid is $128 \times 128$
pixels. Fig.~\ref{acadex} shows three reconstructions based on
\eqref{e-ml-aml} with different values of the parameter $\omega$. The
effect of the flat-field error appears as a ring in the
reconstructions, and it is clear that the severity of both noise and
the ring in the reconstruction decreases as the flat-field intensity
is increased. \remove{The radius of the ring depends on the position of the
detector element for which the flat-field estimate is wrong.} In the
next section, we propose and investigate a new reconstruction model
that takes a statistical model of the flat-field into account.

\begin{figure}
  \centering
  \begin{tikzpicture}[font=\footnotesize, trim axis left, trim axis right]
    \begin{axis}[xmin=-2.0, xmax=2.0, ymin=-5, ymax=10,
      xlabel=$\rho$, ylabel=$\tilde\mu(\rho)$,
      width=0.9\columnwidth,height=0.6\columnwidth,font=\small,xtick={-2,-1.5,-1.0,-0.5,0,0.5,1.0,1.5,2},xticklabels={$-2$,,$-1$,,$0$,,$1$,,$2$},
      legend columns=-1,
      legend style={font=\scriptsize,draw=none,anchor=west,column sep=2mm},
      legend entries={$t_0=0.5$,$t_0=1.0$,$t_0=1.5$},
      legend cell align=left,
      legend to name=named1]
      \addplot[thick,black,smooth] table[x index=0, y index=1]{data/ring_integral.dat};
      \addplot[thick,red,smooth, densely dashed] table[x index=0, y index=2]{data/ring_integral.dat};
      \addplot[thick,blue,smooth, densely dotted] table[x index=0, y index=3]{data/ring_integral.dat};
      \addplot[thin, dotted] table[x index=0, y index=1]{data/ring_envelope.dat};
      \addplot[thin, dotted] table[x index=2, y index=3]{data/ring_envelope.dat};
      \addplot[thin, dotted] table[x index=0, y index=1, x expr=-\thisrowno{0}]{data/ring_envelope.dat};
      \addplot[thin, dotted] table[x index=2, y index=3, x expr=-\thisrowno{2}]{data/ring_envelope.dat};
    \end{axis}
  \end{tikzpicture}
  \ref{named1}
  \caption{Examples of radial profile of reconstruction of
    $h_{\theta}(t)$ for three different values of $t_0$ (0.5, 1.0, and
    1.5) and $\epsilon = 0.05$.}
\label{fig-iradon-sino-delta}
\end{figure}
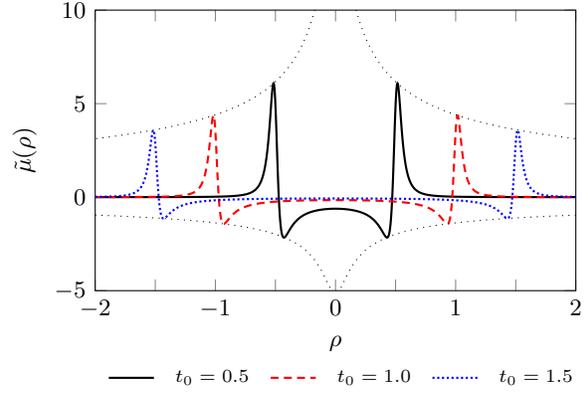

The effect of a flat-field estimation error on the reconstruction may
also be analyzed by means of an analytic reconstruction of the
sinogram $h_{\theta}(t) = \delta(t - t_0)$ where $t_0 \neq 0$ is a
given constant. This corresponds to a ``line'' in the sinogram. The
function $h_{\theta}(t)$ is a radial function (\ie, it does not depend
on $\theta$), but it is not the Radon transform of a function since
$h_{\theta}(t) \neq h_{\theta + \pi}(-t)$. As a consequence, the
Fourier slice theorem does not hold. However, we may still compute a
reconstruction using filtered backprojection. The reconstruction
$\mu(x)$ is itself a radial function, and if we let
$x = \rho n_{\phi}$ where $n_{\phi} = (\cos \phi ,\sin \phi)$ such
that $|\rho|$ is the distance to the origin, we obtain the
expression \cite{KaS:01}
\begin{align*}
  \mu(\rho n_{\phi}) &= \frac{1}{2} \int_{-\pi}^{\pi} \int_{-\infty}^\infty H_{\theta}(\zeta) |\zeta| e^{-2\pi\epsilon
  |\zeta|} e^{i2\pi \zeta \rho n_{\phi}^Tn_\theta }\, d \zeta \,
                       d\theta \\
  &= \pi \int_{-\infty}^\infty H_{\theta}(\zeta) |\zeta| e^{-2\pi\epsilon  |\zeta|}  J_0(2\pi \zeta  \rho) \, d \zeta \\
  &= \pi \int_0^\infty \left[ H_\theta(\zeta) + H_\theta(-\zeta)\right] e^{-2\pi\epsilon \zeta} \zeta
    J_0(2\pi \zeta \rho) \, d\zeta
\end{align*}
where $J_0$ denotes the zeroth-order
Bessel function of the first kind, $H_{\theta}(\zeta) = e^{-i 2\pi \zeta t_0}$ is the Fourier transform of
$h_\theta(t)$, and $|\zeta| e^{-2\pi\epsilon
  |\zeta|}$ is an apodizing filter with parameter $\epsilon
>0$. Using the Hankel transform pair
(20) in \cite[p.~9]{Bat:54}, we obtain the closed-form expression
\begin{align}
  \label{e-delta-fbp}
  \tilde \mu(\rho) &= \frac{1}{4\pi} \left( \frac{\sigma}{(\sigma^2 +
                     \rho^2)^{3/2}} + \frac{\bar\sigma}{(\bar \sigma^2 + \rho^2)^{3/2}} \right) 
\end{align}
where $\sigma = \epsilon + it_0$ and  $\tilde \mu(\rho) = \mu(\rho n_{\phi})$ . 
Fig.~\ref{fig-iradon-sino-delta}  shows three
examples of what this function may look like.
It is clear from the figure that a systematic error in the sinogram in
the form of a ``line'' will appear as spikes in the radial
reconstruction. In particular, the reconstruction will have two
``rings'' of opposite sign near $\rho = t_0$, corresponding to the
positive and negative peaks in the profile $\tilde \mu(\rho)$. The
extrema of $\tilde \mu(\rho)$ (\ie, the spike magnitudes) depend on
both $t_0$ and $\epsilon$. The dotted curves in the figure provide an
envelope of the extrema for $\epsilon = 0.05$, and it shows that the
magnitude of a spike is large when $|t_0|$ is small and vice
versa. Our analysis of the extrema of $\tilde \mu(\rho)$,
\add{which is included in Appendix~\ref{a-extrema},} 
shows that they are approximately
inversely proportional to $\sqrt{\epsilon^3 |t_0|}$ when
$|t_0| \gg \epsilon$. \remove{This is consistent
with the behavior of the rings in Fig.~\ref{acadex} where the
magnitude of each ring depends on its radius.}
Moreover, $\tilde \mu(\rho)$ may have a
significant offset near $\rho = 0$, as is the case for the example
with $t_0=0.5$ in Fig.~\ref{fig-iradon-sino-delta}.

\subsection{Including Prior Information}
If the prior probability density $\prob(u)$ is assumed to be known, a so-called maximum a
posteriori (MAP) estimate can be expressed as
\begin{align}\label{e-map-map}
  \hat u_{\mathrm{map}} = \argmin_{u} \left\{ -\log \prob(u\mid y, v) \right\}
\end{align}
where, according to Bayes' rule, the posterior probability density
$\prob(u| y,v)$ satisfies
\begin{align} 	\label{probpost}
	\prob(u\mid y, v) \propto \prob(y\mid u, v)\prob(u).
\end{align}
Again, since $v$ is generally unknown, an approximate MAP (AMAP) estimate
can be obtained by maximizing an approximation of the posterior
distribution, \ie,
\begin{align}\label{e-map-amap-form}
  \hat u_{\mathrm{amap}} = \argmin_{u} \left\{ -\log \prob(u\mid y,
  \vml) \right\}.
\end{align}
We will restrict our attention to priors of the form
\begin{equation}
	\prob(u\mid\gamma) \propto e^{-\gamma \phi(u)}
		\label{probu}
\end{equation}
where $\phi(u)$ is a convex function and $\gamma > 0$ is a
hyperparameter\remove{, and as a consequence, the estimation problem \eqref{e-map-amap-form} is convex
if the ML estimation problem \eqref{e-ml-aml} is
convex}. \add{With this prior, the AMAP estimation problem can
be expressed as
\begin{align}\label{e-map-amap}
  \hat u_{\mathrm{amap}} = \argmin_{u}  \left\{ (\ones \otimes \vml)^T \exp(-Au) +
    y^TAu + \gamma \phi(u) \right\}
\end{align}
which is a convex optimization problem. Alternatively, using the quadratic
approximation \eqref{e-ml-wls-cost} in place of the log-likelihood
function, we obtain the regularized weighted least-squares problem
\begin{align}\label{e-map-wls}
 \hat u_{\mathrm{wls}} =  \argmin_{u} \left\{ \frac{1}{2} \|
  \diag(y)^{1/2} (Au- b)\|_2^2 + \gamma \phi(u) \right\}
\end{align}
as an approximation to the AMAP estimation problem. }

\section{Joint Reconstruction Approach} \label{proposed} 

We now turn to the main contribution of this paper, namely a model for
jointly estimating the flat-field $v$ as well as the absorption image
$u$. \add{ Recall from the example in section
  \ref{effecffferror} that the approximate ML model
  \eqref{e-ml-aml} may lead to ring artifacts. As will be evident from our
  numerical experiments in section \ref{results}, the approximate MAP model \eqref{e-map-amap}
  suffers the same drawback. To mitigate this, we consider joint MAP
  estimation of $u$ and $v$.} \remove{Instead of solving the
  approximate MAP problem \eqref{e-map-amap}, which we have seen may
  lead to ring artifacts in the reconstruction, we consider the joint
  MAP estimation of $u$ and $v$.} This approach is motivated by the
fact that the measurements $Y$ contain information about both $u$ and
$v$. Indeed, given $u$, an ML estimate of $v$ can be computed as
\begin{align}\label{e-mlu-aml}
 \vmly(u) &= \argmin_{v} \left\{ -\log \prob(Y\mid u, v)
\right\} \\
&= \diag\left( \sum_{j=1}^p \exp(-A_ju) \right)^{-1}Y \ones.
\end{align}

\subsection{MAP Estimation Problem} \label{ss-map-estimation}
With the model assumptions described in \ref{ss-model} and given a
flat-field prior $\prob(v | \alpha, \beta)$, the joint posterior
distribution of the unknown parameters $u$ and $v$ can be expressed as
\begin{align*}
\prob(u,v\mid Y, F) \propto  \prob(Y, F \mid u, v) \prob(u\mid \gamma) \prob(v \mid \alpha,
  \beta)
\end{align*}
where $\prob(Y, F | u, v) = \prob(Y | u,v)\prob(F | v)$, and
$\alpha \in \reals^r$ and $\beta \in \reals^r$ are hyperparameters
associated with the flat-field prior. Here we
will assume that $v_i$ and $v_j$, $i\neq j$ are independent,
and the flat-field prior is $v_i | \alpha_i,\beta_i \sim
\mathrm{Gamma}(\alpha_i,\beta_i)$ for $i=1,\ldots,r$, \ie,
\[ \prob(v_i\mid \alpha_i, \beta_i) =
  \frac{\beta_i^{\alpha_i}}{\Gamma(\alpha_i)} v_i^{\alpha_i -1}
  \exp(-\beta_i v_i). \] The Gamma prior is chosen because of
computational convenience; it is the so-called conjugate prior for the
Poisson likelihood function, and as a consequence, the posterior
distribution of $v$ given $u$ is itself a Gamma distribution. For the
Gamma distribution, the hyperparameter $\alpha_i$ is commonly referred
to as the shape, and $\beta_i$ is referred to as the rate.  The
\remove{resulting}\add{ corresponding} MAP estimation problem can be expressed as
\begin{align}\label{e-map-jmap}
  (\hat u,\hat v) &= \argmin_{(u,v)} \left\{ -\log \prob(u,v\mid Y, F)
  \right\} \\
\notag   &= \argmin_{(u,v)} \left\{ J(u,v) + \gamma \phi(u)\right\}
\end{align}
where
\begin{align}
  J(u,v) = v^Td(u)  + y^TAu -c^T\log(v)
\end{align}
and
\begin{align}
  \label{e-c-du}
 c = F \ones + Y\ones + \alpha - \ones,\  d(u) = s\ones + \sum_{j=1}^p \exp(-A_ju) + \beta.  
\end{align}
The function $J(u,v)$ is convex in $u$ given $v$ and vice versa, but
it is not jointly convex in $u$ and $v$. However, by setting the gradient
of $J(u,v)$ with respect to $v$ equal to zero, we obtain the
first-order optimality condition
$ \hat v(u) = \diag(d(u))^{-1} c$.
This allows us to eliminate $v$ from the estimation problem
\eqref{e-map-jmap}, \ie,
\begin{align*}
J(u,\hat v(u)) \propto  y^TAu + c^T\log(d(u)),
\end{align*}
which is a convex function of $u$. \remove{Thus, the problem
\eqref{e-map-jmap} can be solved by computing the MAP estimates}
\add{ Thus, the problem
\eqref{e-map-jmap} is equivalent to the following convex reconstruction model}
\begin{align} \label{e-map-umap}
 \hat u &=  \argmin_{u} \left\{   y^TAu + c^T\log(d(u)) +  \gamma
          \phi(u) \right\}
\end{align}
\add{ with the flat-field estimate $\hat v$ given by}
\begin{align}
\label{e-map-vmap}
\hat v &=  \diag(d(\hat u))^{-1} c.
\end{align}
We note that \remove{the flat-field estimate} $\hat v$ has an interesting interpretation:
each element of \remove{the flat-field estimate} $\hat v$ can be expressed as a
convex combination of three independent estimates, \ie,
\begin{align} \label{theta-par}
  \hat v = \diag(\theta_1) \vml + \diag(\theta_{2}) \vmly(\hat u) + \diag(\theta_3) \hat v_{\mathrm{pr}}(\alpha,\beta)
\end{align}
where $\theta_1,\theta_2,\theta_3 \in \reals_+^{r}$,
$\theta_1+\theta_2+\theta_3 = \ones$, are parameters that depend on both data
and $\hat u$, $\alpha$, and $\beta$. The ML estimate $\vml$, defined
in \eqref{empmeanest}, is based on the flat-field measurements $F$,
the estimate $\vmly(\hat u)$ is based on the measurements $Y$ and
defined in \eqref{e-mlu-aml}, and the estimate $\hat
v_{\mathrm{pr}}(\alpha,\beta) = \diag(\beta)^{-1}(\alpha - \ones)$ is
based on the flat-field prior\add{; see Appendix \ref{flat-field-intr} for further details on this interpretation.}

\subsection{Choosing The Hyperparameters}\label{ss-hyperparameters}

The estimation problem \eqref{e-map-umap} depends on the flat-field
hyperparameters $\alpha$ and $\beta$. We now discuss different ways to
choose these hyperparameters.

\subsubsection{Uniform Positive Prior}
The simplest prior is perhaps the uniform positive (UP) prior which is
obtained by setting $\alpha_i = 1$ and $\beta_i = 0$ for
$i=1,\ldots,r$. In the present case, this corresponds to simply omitting the prior
$\prob(v|\alpha,\beta)$ from the model, and hence the flat-field estimates $\hat
v(u)$ become convex combinations of only two estimates instead of
three. This is an improper prior since it does not integrate to one.

\subsubsection{Jeffreys Prior}
The Jeffreys prior (JP) for the Poisson distribution is
$p( v_i | \alpha_i,\beta_i) \propto 1/\sqrt{v_i}$ which is obtained
by letting $\alpha_i= 0.5$ and $\beta_i=0$. This is also an improper
prior.

\subsubsection{Type-II ML Estimation}

The flat-field measurements can be used to estimate the
hyperparameters by \remove{computing and} maximizing the marginal probability
of $f_{i1},\ldots,f_{is}$ given the hyperparameters $\alpha_i$ and
$\beta_i$, \ie,
\begin{align}\label{e-ml-typeII}
  (\hat \alpha_i, \hat\beta_i) = \add{\argmin_{(\alpha_i,\beta_i)}}\remove{\argmin_{(\alpha,\beta)}} \left\{-\log \prob(f_{i1},\ldots,f_{is} \mid \alpha_i, \beta_i) \right\}.
\end{align}
This is known as type-II ML estimation or empirical Bayes estimation
\cite{Ber:85}. \add{As shown in Appendix
  \ref{type-II-estimates}, this approach leads to the AMAP model, \ie,
a zero-variance prior with mean $\vml$.}

\subsubsection{Flat-field Emphasizing Prior}  \label{EmpMean}
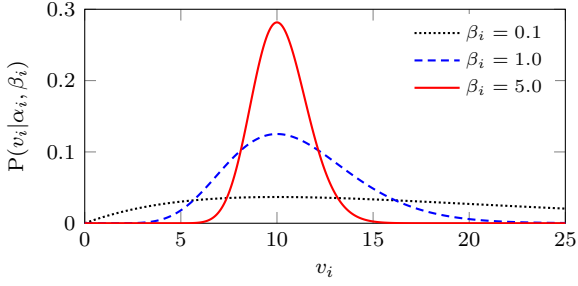
\begin{figure}
  \centering
  \begin{tikzpicture}[font=\footnotesize, trim axis left, trim axis right]
    \begin{axis}[no markers,xmin=0,xmax=25,ymin=0,ymax=3e-1, width=0.9\columnwidth,height=0.5\columnwidth,
      ylabel={$\prob(v_i|\alpha_i,\beta_i)$},xlabel={$v_i$},grid=none,legend
      style={font=\scriptsize, draw=none,anchor=north east},
      legend entries={$\beta_i = 0.1$,$\beta_i=1.0$,$\beta_i=5.0$},
      legend cell align=left]
      \addplot[thick, black, smooth, densely dotted] table [x index=0, y index=1]{data/gamma_examples.dat};
      \addplot[thick, blue,smooth,densely dashed] table [x index=0, y index=2]{data/gamma_examples.dat};
      \addplot[thick, red,smooth] table [x index=0, y index=3]{data/gamma_examples.dat};
    \end{axis}
  \end{tikzpicture}
  \caption{Gamma distributions with hyperparameters $\beta_i$ and
    $\alpha_i = 1+ (\vml)_i\beta_i$ for $(\vml)_i = 10$ and $\beta_i \in  \{0.1,1.0,5.0\}$.}
  \label{fig-gamma-examples}
\end{figure}

Recall that the flat-field estimate $\hat v(u)$ can be expressed as
convex combinations of three estimates. Specifically,
\begin{align}\label{e-map-vmap-cc3}
 \hat  v_i(u) &= \frac{s}{d_i(u)} (\vml)_i + \frac{\tau_i(u)}{d_i(u)} (\vmly)_i +
  \frac{\beta_i}{d_i(u)} \frac{\alpha_i-1}{\beta_i}
\end{align}
where $\tau_i(u) = \sum_{j=1}^p \exp(-e_i^TA_ju)$. If we set the
mode of the Gamma prior (\ie, $(\alpha_i - 1)/\beta_i$) equal to the
flat-field ML estimate $(\vml)_i$ by letting
$\alpha_i = 1+\beta_i (\vml)_i$, we obtain the estimate
\begin{align}
  \label{e-map-vmap-cc2}
  \hat  v_i(u) = \frac{s + \beta_i}{d_i(u)} (\vml)_i + \frac{\tau_i(u)}{d_i(u)} (\vmly)_i
\end{align}
which is a convex combination of two estimates. It is easy to verify
that $\hat v_i(u) \rightarrow (\vml)_i$ as $\beta \rightarrow \infty$,
and with $\beta_i = 0$, the \add{estimate $\hat
  v_i(u)$}\remove{FE prior} is equivalent to the \add{
  estimate obtained with the }UP prior.
Thus, choosing $\beta_i > 0$ and $\alpha_i = 1+\beta_i (\vml)_i$
allows us to emphasize the flat-field ML estimate $(\vml)_i$. This is
consistent with the fact that the parameter $\beta_i$ is the rate
parameter associated with the Gamma distribution: the larger the rate,
the more concentrated the distibution is around its mode. This is
illustrated in Fig.~\ref{fig-gamma-examples}. We call this
corresponding prior the flat-field emphasizing (FE) prior.

\subsection{Quadratic Approximation}
A quadratic approximation of the first two terms in \eqref{e-map-umap}
can be derived by means of a second-order Taylor expansion with
respect to $Au$. Substituting $y$ for
$(I\otimes \diag(\hat v(u)))\exp(-Au)$, we obtain the following
approximate MAP estimation problem
\begin{align}\label{e-amap-swls}
  \hat u_{\mathrm{swls}} = \argmin_{u} \left\{ \frac{1}{2}\|Au -
  b\|_{\widehat \Sigma_{\rv{b}}^{-1}}^2 + \gamma \phi(u) \right\}
\end{align}
where the covariance matrix $\widehat \Sigma_{\rv{b}}$ is defined as
\begin{align}
\widehat \Sigma_{\rv{b}} = ( \ones\ones^T) \otimes \diag(s\vml + \alpha -1)^{-1} + \diag(y)^{-1}.
\end{align}
This is also the covariance matrix associated with $\rv{b}$ in the
linear approximation \eqref{e-linear-approximation}.
Note that the weighted least-squares data fidelity term takes
the systematic errors induced by flat-field estimation errors into
account without explicitly modeling the flat-field, and hence we label
this a regularized stripe-weighted least-squares (SWLS) problem. The
model depends on the hyperparameter vector $\alpha$, which appears in the
covariance matrix, but $\beta$ does not appear in the model. 


\section{Implementation}\label{s-implementation}
The MAP estimation problems \eqref{e-map-amap} and
\eqref{e-map-umap}\remove{,} \add{as well as}
the \add{WLS \eqref{e-map-wls} and } SWLS \eqref{e-amap-swls}  \add{ quadratic} approximation\add{s} are all convex problems
that can be solved with a wide range of numerical optimization
methods. Here we will focus on \add{simple} first-order methods which are suitable
for large-scale problems.

\subsection{Attenuation Priors}
Before we describe our implementation of the different reconstruction
methods, we briefly discuss two attenuation priors of the form
\eqref{probu}, namely the nonnegativity prior (corresponding to
nonnegativity constraints $u_i \geq 0$), and a combination of the
nonnegativity prior and total variation (TV) regularization
\cite{Rudin1992}. Both of these priors can be combined with the
existing AMAP model \eqref{e-map-amap}, the proposed model
\eqref{e-map-umap}, \add{the WLS model \eqref{e-map-wls}} and the SWLS model \eqref{e-amap-swls}.

\subsubsection{Nonnegativity}
The nonnegativity constraints can be expressed as $\phi(u) = I_+(u)$
where $I_+(u)$ denotes the indicator function of the nonnegative
orthant, \ie, $I_+(u) = 0$ if and only if $u$ is a nonnegative vector,
and otherwise $I_+(u) = \infty$.

\subsubsection{Nonnegativity and TV}
The combination of nonnegativity constraints and TV may be expressed as
\[ \phi(u) = I_+(u) +  \mathrm{TV}_{\delta}(u) \]
where
$ \mathrm{TV}_{\delta}(u) = \sum_{i=1}^{n} \xi_{\delta}
(\|D_iu\|_2) $
is a differentiable TV-approximation, $\xi_{\delta}$
denotes the Huber-norm
\[ \remove{\xi(t)} \add{\xi_{\delta}(t)}=
  \begin{cases}
    (t)^2/(2\delta)      & \quad |t| \leq \delta\\
    |t| - \delta/2  & \text{otherwise}\\
  \end{cases}
  \]
with parameter $\delta$, and $D_iu$ is a finite-difference
approximation of the gradient at pixel $i$. We will use a pixel basis
corresponding to an $M\times N$ grid (\ie, $n = MN$). Specifically, we
define 
\[ D_i =
  \begin{bmatrix}
    e_i^T(I_N \otimes \bar D_M) \\
    e_i^T(\bar D_N \otimes I_M) 
  \end{bmatrix}\] where $I_M$ and $I_N$ are identity matrices, and $\bar D_M$ and $\bar D_N$ are square
difference matrices of order $M$ and $N$, respectively, and of the
form
\[
  \begin{bmatrix}
      1 & -1 &             &\\
         & \ddots & \ddots & \\
         &              &  1         & -1 \\
         &               &             &  0 
  \end{bmatrix}
\]
where the last row is zero, corresponding to Neumann boundary conditions.

 The function $\mathrm{TV}_{\delta}(u)$ has a
Lipschitz continuous gradient with constant $L_{\mathrm{tv}}(\delta)
= \|D\|_2^2/\delta$ where $D =
\begin{bmatrix}
  D_1^T & \cdots D_n^T
\end{bmatrix}^T$.

\subsection{Reconstruction Models}
We now consider \remove{four} \add{five} different reconstruction models of the form
\begin{align}
  \label{e-rec-models}
  \begin{array}{ll}
    \mbox{minimize} & J_i(u)  + \gamma \phi(u), \quad i =1,\ldots,\remove{4}\add{5},
  \end{array}
\end{align}
where $J_i(u)$ is based on either \eqref{e-map-map}, \eqref{e-map-amap},
\eqref{e-map-umap}, \add{\eqref{e-map-wls}} or \eqref{e-amap-swls}.

\subsubsection{Baseline and AMAP Estimation}
The reconstruction model \eqref{e-map-map} requires the true
flat-field $v$ which is not available in practice. However, the model
may be used to compute a baseline reconstruction in simulation
studies. The baseline reconstruction problem corresponds to
$J_1(u) = J(u,v)$ where the true flat-field $v$ is assumed to be
known. If we replace $v$ by $\vml$, we obtain the
AMAP model \eqref{e-map-amap} with objective
$J_2(u) = J(u, \vml)$.

To solve the reconstruction problem \eqref{e-rec-models} using a
first-order method, we need the gradient of $J(u,v)$ with respect to
$u$, \ie,
\begin{align}\label{e-map-amap-grad}
  \nabla_u  J(u,v) = A^T(y-\hat y(u,v))
\end{align}
where $\hat y(u,v) = (I\otimes \diag(v)) \exp(-Au)$.
It is easy to verify that the gradient $\nabla_u  J(u,v)$ is Lipschitz continuous on the
nonnegative orthant since the norm of the Hessian
\begin{align*}
    \nabla_u^2  J(u,v) = A^T\diag(\hat y(u,v)) A
\end{align*}
is bounded for $u \geq 0$ and with $v$ fixed. We will use the Lipschitz constants $L_1 = \max_i\{v_i\}\|A\|_2^2$ and $L_2 =
\max_i\{(\vml)_i\} \|A\|_2^2$.

\subsubsection{Joint MAP Estimation}
The MAP estimation problem \eqref{e-map-umap} is a special case of
\eqref{e-rec-models} if we let $J_3(u) = J(u, \hat v(u))$. The
gradient of $J_3(u)$ is
\begin{align}\label{e-map-umap-grad}
  \nabla J_3(u) &= A^Ty + D_d(u)^T\hat v(u) \\
  \notag
  &= A^T( y-\hat y(u,\hat v(u)) )
\end{align}
where $D_d(u) = -\sum_{j=1}^p \diag(\exp(-A_ju)) A_j$ denotes the
Jacobian matrix of $d(u)$. Comparing with \eqref{e-map-amap-grad}, we
see that the only difference is that the residual
$y-\hat y(u,\hat v(u))$ is based on the flat-field estimate
$\hat v(u)$ instead of the true flat-field $v$ or the ML estimate
$\vml$.

To derive the Hessian of $J_3(u)$, note that
\[ c^T\log(d(u)) = \sum_{i=1}^r c_i \log(d_i(u)) \]
where $d_i(u) = s + \sum_{j=1}^p \exp(-e_i^TA_ju) + \beta_i$. This
implies that the Hessian can be expressed as
\[ \sum_{i=1}^r c_i \left(
  \frac{\nabla^2 d_i(u) }{d_i(u)}- \frac{ \nabla d_i(u)
  \nabla d_i(u)^T}{d_i(u)^2} \right). \]
Now let $\Pi_i = I \otimes e_i^T$ such that $\Pi_i y = Y^T e_i$ corresponds to
the $i$th row of $Y$, and define a permutation matrix $\Pi = [\Pi_1^T
\ \cdots\ \Pi_r^T ]^T$.
This allows us to express the Hessian $\nabla^2 J_3(u)$ as
\begin{align}\label{e-j3-hessian}
  \nabla^2 J_3(u) = A^T\Pi^T\blkdiag(B_1(u),\ldots,B_r(u))\Pi A
\end{align}
where
$ B_i(u) = \diag(\Pi_i \hat y) - \frac{1}{c_i} \Pi_i \hat
y\hat y^T \Pi_i^T$,
and where $\hat y$ is used as shorthand for $\hat y(u,\hat v(u)
)$. (We remark that $\hat v(u)$ depends on both $\alpha$ and
$\beta$, and consequently, so does the Hessian $\nabla^2 J_3(u)$.)  It follows that
\[\|\nabla^2 J_3(u)\|_2 \leq \|A^T\diag(y)A\|_2 \]
which implies that $\nabla J_3(u)$ is Lipschitz continuous with
constant $L_3 = \|A^T\diag(y)A\|_2$.

\add{
\subsubsection{WLS Estimation}
The quadratic approximation \eqref{e-map-wls} corresponds to \eqref{e-rec-models} with 
$J_4(u) = \frac{1}{2} \| Au - b\|_{\widehat \Sigma_{\rv{b}}^{-1}}^{2}$
and
$ \widehat \Sigma_{\rv{b}} = \diag(y)^{-1} $.
The gradient of $J_4(u)$ is
\[ \nabla J_4(u) = A^T\Sigma_{\rv{b}}^{-1}(Au-b) \]
which is Lipschitz continuous with constant
$\|A^T\widehat\Sigma_{\rv{b}}^{-1}A \|_2$.
}

\subsubsection{Regularized SWLS}
The quadratic approximation \eqref{e-amap-swls} corresponds to
\eqref{e-rec-models} with
$\remove{J_4(u)}\add{J_5(u)} = \frac{1}{2} \| Au - b\|_{\widehat \Sigma_{\rv{b}}^{-1}}^{2}$
and
\[ \widehat \Sigma_{\rv{b}} =\Pi^T \left[ \diag(\Pi y)^{-1}  + \diag(s\vml + \alpha -\ones)^{-1}
\otimes (\ones \ones^T) \right]\Pi.  \]
Thus, $\widehat \Sigma_{\rv{b}}$ is a symmetric permutation
of a block-diagonal matrix with diagonal-plus-rank-one blocks, and hence
matrix-vector products with $\widehat \Sigma_{\rv{b}}^{-1} $ can be
efficiently evaluated using the Woodbury identity, \ie,
$
  \widehat \Sigma_{\rv{b}}^{-1} = \Pi^T \blkdiag(S_1,\ldots,S_r)\Pi
$
where
\begin{align}\label{e-swls-blocks}
  S_i = \diag(\Pi_iy) - \frac{1}{s(\vml)_i + e_i^{T}Y\ones + \alpha_i -1} \Pi_iyy^T\Pi_i^T.
\end{align}
This allows us to evaluate the gradient as
\[ \nabla \remove{J_4(u)}\add{J_5(u)} = A^T\Sigma_{\rv{b}}^{-1}(Au-b) \]
which is Lipschitz continuous with constant
$\|A^T\widehat\Sigma_{\rv{b}}^{-1}A \|_2$.

It is instructive to compare the SWLS \remove{mode}\add{model} to the WLS model considered in
\cite{AdityaMohan2015}. This model implicitly includes the
flat-fields using the following objective function
\begin{align}\label{e-wls-z}
\remove{J_5(u,z)}\add{J_6(u,z)}= \frac{1}{2} \|
  \diag(y)^{1/2} (Au- b + \ones\otimes z)\|_2^2 + \frac{\lambda}{2}\| z \|_2^2
\end{align}
where $z \in \reals^r$ is an auxiliary variable that can be thought of
as the relative flat-field error (cf.\ the analysis in Section
\ref{effecffferror}).  Taking the gradient
with respect to $z$ and setting it equal to zero yields
$
  z =  \diag(Y\ones + \lambda \ones)^{-1}(\ones^T \otimes I)\diag(y) (b-Au),
$
and using this expression in \eqref{e-wls-z} yields
\begin{align}\label{e-wls-zelim}
\remove{J_5(u)}\add{J_6(u)}=  \frac{1}{2} \| Au - b\|_{\widehat \Sigma^{-1}}^{2}
\end{align}
where $ \widehat \Sigma^{-1} = \Pi^T \blkdiag(\bar S_1,\ldots,\bar S_r)\Pi$
and
\begin{align}\label{e-wls-blocks}
  \bar S_i = \diag(\Pi_iy) - \frac{1}{e_i^{T}Y\ones + \lambda} \Pi_iyy^T\Pi_i^T.
\end{align}
The blocks $\bar S_i$ clearly resemble the blocks $S_i$ from the SWLS
model in \eqref{e-swls-blocks}: the only difference is the scalar
weight in front of the rank-1 term in each of the $r$ blocks. In
particular, notice that the weights in the SWLS model include
information derived from all measurements as well as the flat-field
prior. Moreover, the parameter $\lambda$ in \eqref{e-wls-blocks} plays
a similar role as the flat-field hyperparameters $\alpha$ in
\eqref{e-swls-blocks}, but the SWLS model is more general and flexible
because it allows the use of a different hyperparameter $\alpha_i$ for
each of the $r$ blocks\remove{, accomodating the type-II ML estimation of Sec.\ \ref{ss-hyperparameters}}.

\subsection{Algorithm}\label{optmethod}
The functions $J_1(u),\ldots,\remove{J_4(u)}\add{J_5(u)}$ are all differentiable with
Lipschitz continuous gradients on the nonnegative orthant, and hence
we can apply a proximal gradient method which is suitable for
minimizing problems of the form
\begin{align*}
  \begin{array}{ll}
    \mbox{minimize} & g(u) + h(u).
  \end{array}
\end{align*}
Here $g : \reals^n \to \reals$ is convex with a Lipschitz continuous
gradient with Lipschitz constant $L$, $h : \reals^n \to \reals$ is convex,
and the prox-operator
\[ \prox_{t h}(\bar u) = \argmin_{u} \left\{ t h(u) + \frac{1}{2} \|u
  - \bar u\|_2^2 \right\} \]
is assumed to be cheap to evaluate. We will define
$g(u) = J_i(u) + \gamma \mathrm{TV}_{\delta}(u)$ and
$h(u) = I_+(u)$, and hence the Lipschitz constant is given by
$L = L_i + \gamma L_{\mathrm{tv}}(\delta)$. Given a starting point
$u^{(0)}$ and a fixed number of iterations $K$, the algorithm can be
summarized as
\[ u^{(k)} = \prox_{t h} ( u^{(k-1)} - t \nabla g (u^{(k-1)})
), \ k=1,2,\ldots,K \]
where \remove{$t=1/L$} \add{$t \in (0,2/L)$} is the step size and
$ \prox_{t h} ( \bar u ) = \max(0,\bar u) $
is the projection onto the nonnegative orthant. With this step size,
the method is a descent method. The Lipschitz constant $L$ can be
estimated without an explicit representation of $A$ or $D$ by means of
the power iteration algorithm. \add{Our MATLAB implementation of the method is available for download at \url{https://github.com/hariagr/R2CT}.}

\section{Numerical Experiments}  \label{results}

\subsection{Simulation Study}
To evaluate the proposed reconstruction models, we conducted a series
of experiments in MATLAB based on simulated data. In these
experiments, we used a parallel beam geometry with $p=720$ equidistant
projection angles covering half a rotation, and a 2 cm wide photon
counting detector array with $r=512$ detector elements. To model a
non-uniform detection efficiency, the elements
of the flat-field vector $v$ were drawn from a Poisson distribution
with mean $I_0$. We used $s = 5$ measurements of the flat-field which
were generated according to \eqref{e-meas-model-ff}, and the
measurements $Y$ were generated according to \eqref{e-meas-model}
using a $2N\times 2N$ pixel discretization of a 2D phantom defined on
a 4 cm$^2$ square. To avoid inverse crimes, we computed our
reconstructions on an $N \times N$ \add{($N = 512$)} pixel grid with a circular
mask. The value of the TV-smoothing parameter $\delta$ was 0.01
cm$^{-1}$ in
all experiments with the TV-prior. \add{We used as step size
  $t=1.8/L$, and we} \remove{We} used the ASTRA Toolbox \cite{APB+:15}
(version 1.7.1beta) to compute filtered backprojection (FBP)
reconstructions and to implicitly compute products with $A$ and $A^T$
on a GPU. We generated the phantoms using the AIR Tools package
\cite{Hansen2012} (version 1.3), and we used the method outlined in Section
\ref{optmethod} to numerically solve the reconstruction problems.  As
a remark, we note that the ASTRA GPU code for backprojection (\ie,
multiplication by $A^T$) is not an exact adjoint of the forward
operator (multiplication by $A$), and this may introduce small errors
in the gradient computations. However, it is significantly faster than
matched implementations, and we did not see any noticeable
differences in reconstruction quality when using the exact
adjoint. 

As initial guess we used a vector of zeros, and we used a fixed number
of iterations as stopping criteria (500 iterations for reconstructions
without the TV-prior and 1,500 iterations for reconstructions with the
TV-prior).  We determine the parameter $\gamma$ for the TV-prior based
on the subjective visualization. As flat-field prior
$\prob(v|\alpha,\beta)$ we used $\alpha_i = 1 + \beta_i (\vml)_i$ and
$\beta_i \geq 0$ (corresponding to the UP flat-field prior if
$\beta_i=0$ and the FE prior if $\beta_i >0$), and for the attenuation
prior $\prob(u|\gamma)$ we used either nonnegativity or total variation
combined with nonnegativity. Note that SWLS only depends on $\alpha$,
but since we also use $\alpha_i = 1 + \beta_i (\vml)_i$ for SWLS, we
report the value of $\beta$ in the experiments.

To quantitatively compare the quality of reconstructions, we
report the relative attenuation error (RAE)
\[
 e_{u}^{\mathrm{rel}}(\hat u) = 100\cdot \frac{\| \hat u -
  u\|_2}{\|u\|_2}, \]
the relative flat-field error (RFE)
\[
   e_{v}^{\mathrm{rel}}(\hat v) = 100\cdot \frac{\|\hat v - v\|_2}{\|v\|_2},
\]
the structural similarity (SSIM) index\footnote{We used the MATLAB
  \texttt{ssim} function with the radius parameter equal to 0.2 for
  reconstructions without the TV-prior and equal to 2.0 for reconstructions with the TV-prior.}
\cite{WBSS:04}, and a ``ring ratio'' (RR), defined as  \[\|\psi_v(\hat
  v(\hat u))\|_F/\|\psi_v(\vml)\|_F\] with $\psi_v(\hat v)$ defined as
\begin{align}
  \label{e-ring-err}
  \psi_v(\hat v) = \mathrm{FBP}( \diag(v)^{-1}(\hat v
  - v) \ones^T )
\end{align}
and where FBP denotes the filtered backprojection reconstruction method.
In other words, $\psi_v(\hat v)$ is the FBP reconstruction of the sinogram
stripes due to flat-field estimation errors, and hence the norm
$\|\psi_v(\hat v)\|_F$ quantifies how severely the flat-field estimation errors affect the
reconstruction. Thus, the RR can be viewed as an indication of the expected ring
artifact reduction if we were to use the flat-field estimate
$\hat v(\hat u)$ instead of the ML estimate $\vml$ (smaller is better)
to compute a reconstruction. Recall that all but the JMAP
reconstruction model are based on the ML estimate $\vml$, so for the
other models, the RFE and the RR simply
reflect what we obtain if we were to use the reconstruction
$\hat u$ to compute a new flat-field estimate $\hat v(\hat u)$\add{, using \eqref{e-map-vmap}. We used $\alpha = 1$ and $\beta = 0$ to compute $\hat v(\hat u)$ for all but the JMAP and SWLS reconstruction models}.

\subsubsection{Low Intensity}
\def\magn{2}
\begin{figure*}
\centering

  \begin{tikzpicture}[x=0.33\columnwidth,y=0.33\columnwidth,font=\footnotesize,spy scope={magnification=\magn, size=0.12\columnwidth,every spy in node/.style={draw}}]

  \node[label={[font=\footnotesize]above:Phantom},inner sep=0pt] (a1) at (0,0) {\includegraphics[width=0.32\columnwidth]{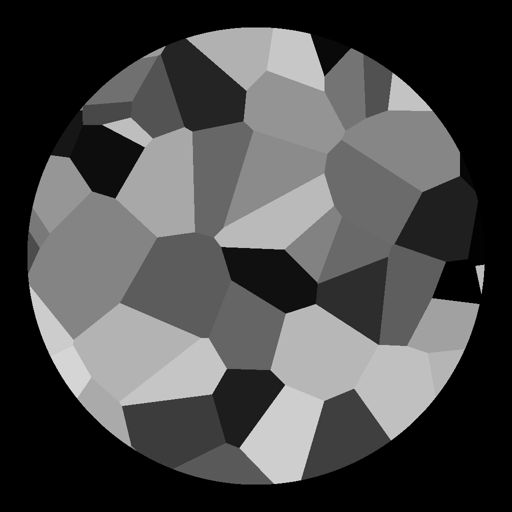}};
  \spy on (a1) in node at (0.3,-0.3);

  \node[label={[font=\footnotesize]above:Baseline FBP},inner sep=0pt] (a2) at (1,0) {\includegraphics[width=0.32\columnwidth]{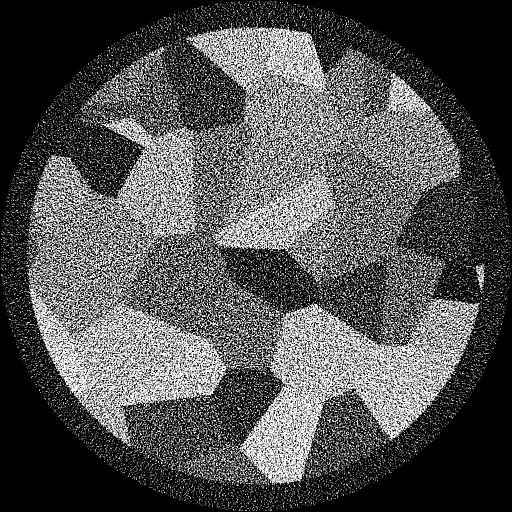}};
  \spy on (a2) in node at (1.3,-0.3);

  \node[label={[font=\footnotesize]above:FBP},inner sep=0pt] (a3) at (2,0) {\includegraphics[width=0.32\columnwidth]{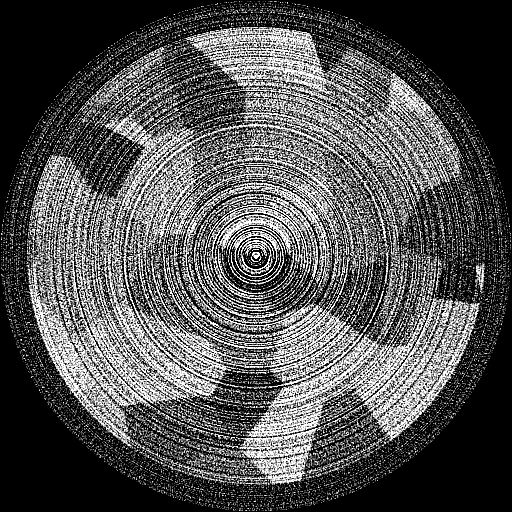}};
  \spy on (a3) in node at (2.3,-0.3);

  \node[label={[font=\footnotesize]above:P-FBP},inner sep=0pt] (a4) at (3,0) {\includegraphics[width=0.32\columnwidth]{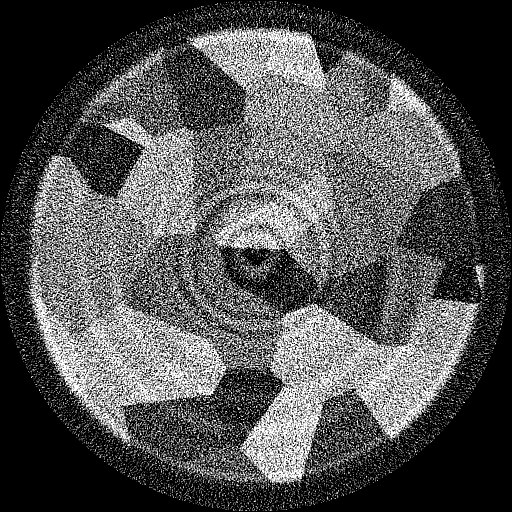}};
  \spy on (a4) in node at (3.3,-0.3);

  \node[label={[font=\footnotesize]above:WLS},inner sep=0pt] (a5) at (4,0) {\includegraphics[width=0.32\columnwidth]{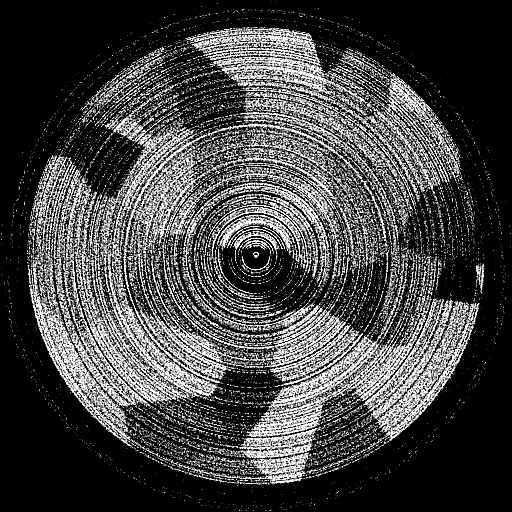}};
  \spy on (a5) in node at (4.3,-0.3);

  \node[label={[font=\footnotesize]above:WLS-TV},inner sep=0pt] (a6) at (5,0) {\includegraphics[width=0.32\columnwidth]{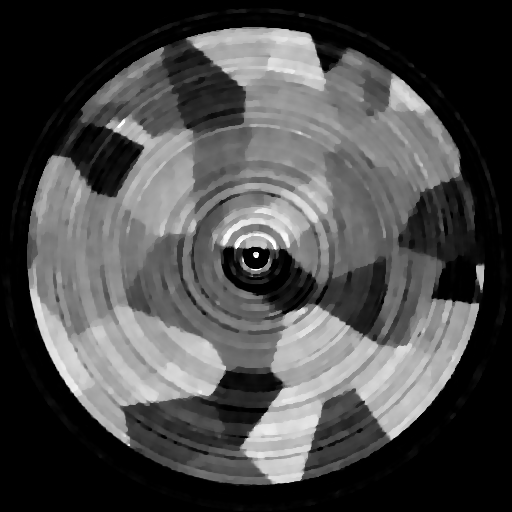}};
  \spy on (a6) in node at (5.3,-0.3);

  \node[label={[font=\footnotesize,rotate=90,anchor=south]left:Without
    TV prior},label={[font=\footnotesize]above:Baseline MAP},inner sep=0pt] (b1) at (0,-1.2) {\includegraphics[width=0.32\columnwidth]{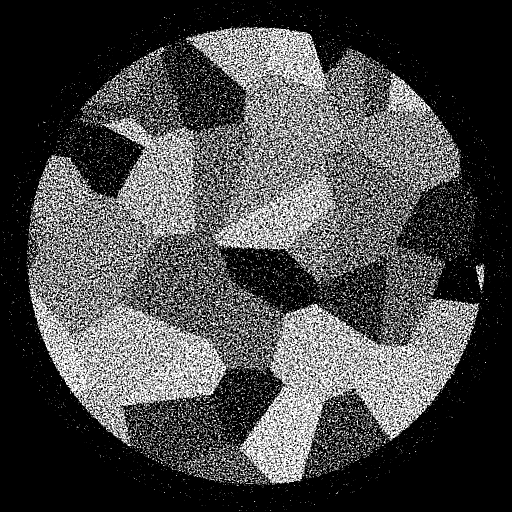}};
  \spy on (b1) in node at (0.3,-1.5);

  \node[label={[font=\footnotesize]above:AMAP},inner sep=0pt] (b2) at (1,-1.2) {\includegraphics[width=0.32\columnwidth]{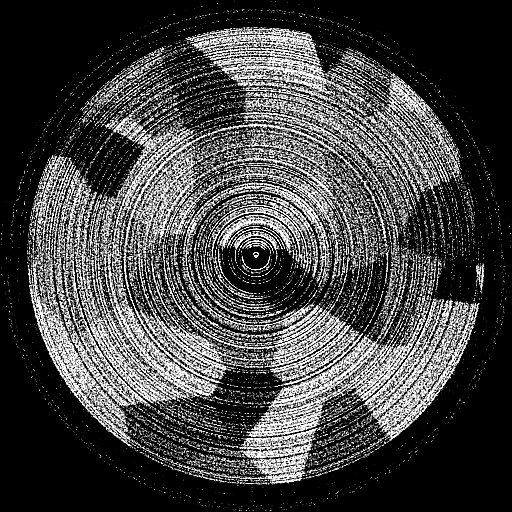}};
  \spy on (b2) in node at (1.3,-1.5);

  \node[label={[font=\footnotesize]above:JMAP  ($\beta=0$)},inner sep=0pt] (b3) at (2,-1.2) {\includegraphics[width=0.32\columnwidth]{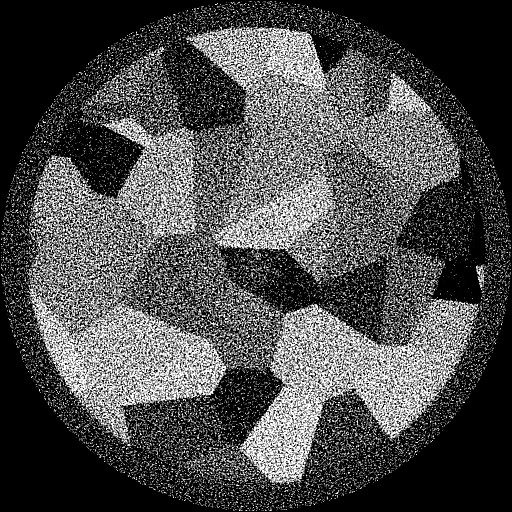}};
  \spy on (b3) in node at (2.3,-1.5);

  \node[label={[font=\footnotesize]above:SWLS  ($\beta=0$)},inner sep=0pt] (b4) at (3,-1.2) {\includegraphics[width=0.32\columnwidth]{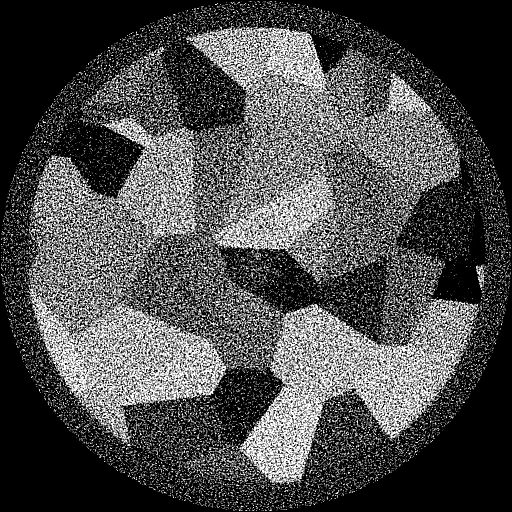}};
  \spy on (b4) in node at (3.3,-1.5);

   \node[label={[font=\footnotesize]above:JMAP  ($\beta=10$)},inner sep=0pt] (b5) at (4,-1.2) {\includegraphics[width=0.32\columnwidth]{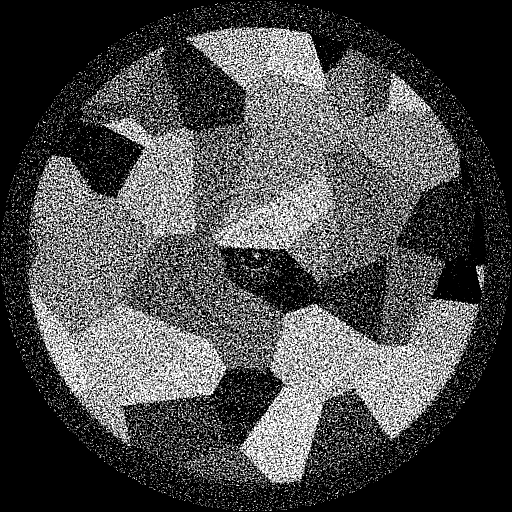}};
  \spy on (b5) in node at (4.3,-1.5);

  \node[label={[font=\footnotesize]above:SWLS  ($\beta=10$)},inner sep=0pt] (b6) at (5,-1.2) {\includegraphics[width=0.32\columnwidth]{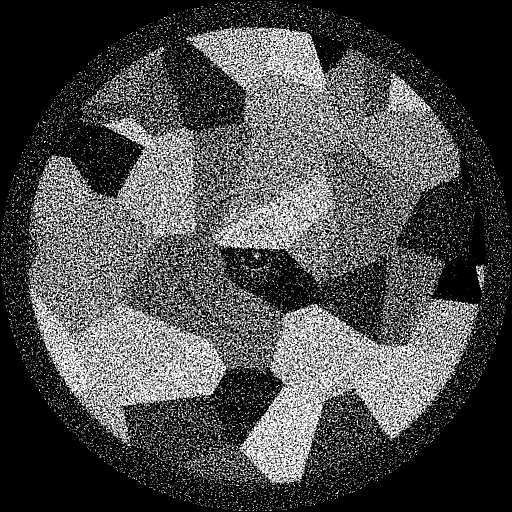}};
  \spy on (b6) in node at (5.3,-1.5);

  \node[label={[font=\footnotesize,rotate=90,anchor=south]left:With TV
  prior},inner sep=0pt] (c1) at (0,-2.2) {\includegraphics[width=0.32\columnwidth]{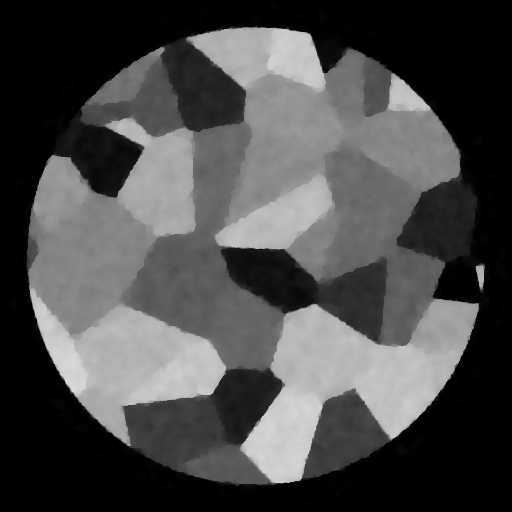}};
  \spy on (c1) in node at (0.3,-2.5);

  \node[label,inner sep=0pt] (c2) at (1,-2.2) {\includegraphics[width=0.32\columnwidth]{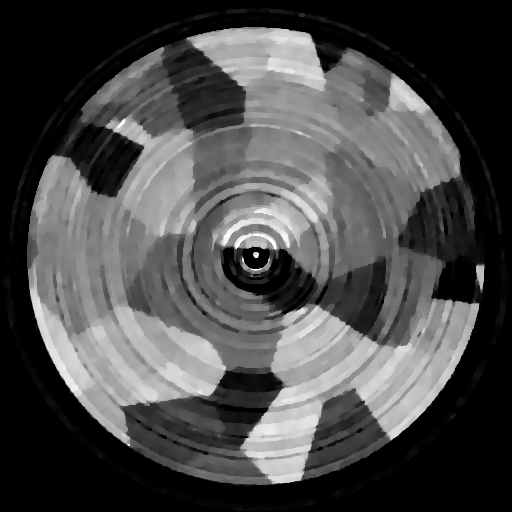}};
  \spy on (c2) in node at (1.3,-2.5);

  \node[label,inner sep=0pt] (c3) at (2,-2.2) {\includegraphics[width=0.32\columnwidth]{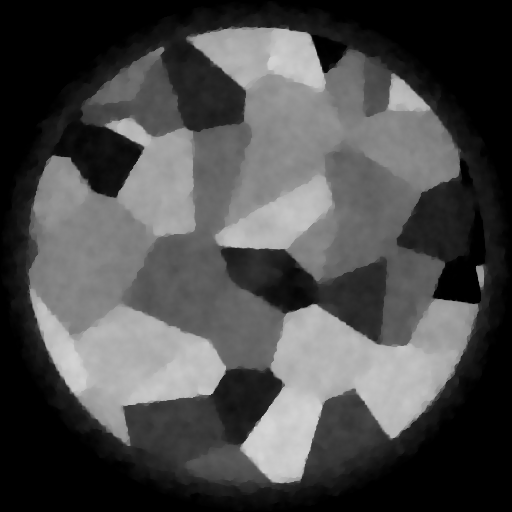}};
  \spy on (c3) in node at (2.3,-2.5);

  \node[label,inner sep=0pt] (c4) at (3,-2.2) {\includegraphics[width=0.32\columnwidth]{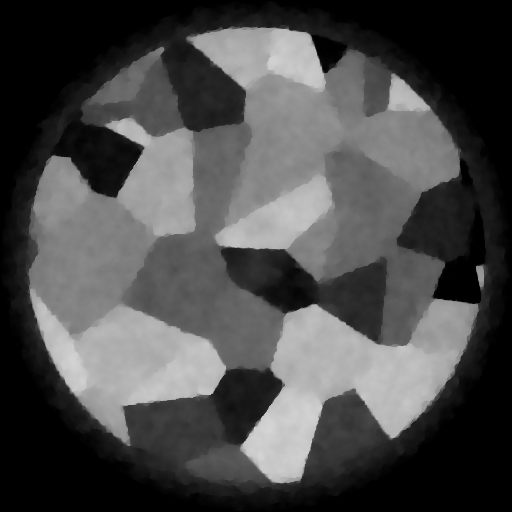}};
  \spy on (c4) in node at (3.3,-2.5);

\node[label,inner sep=0pt] (c5) at (4,-2.2) {\includegraphics[width=0.32\columnwidth]{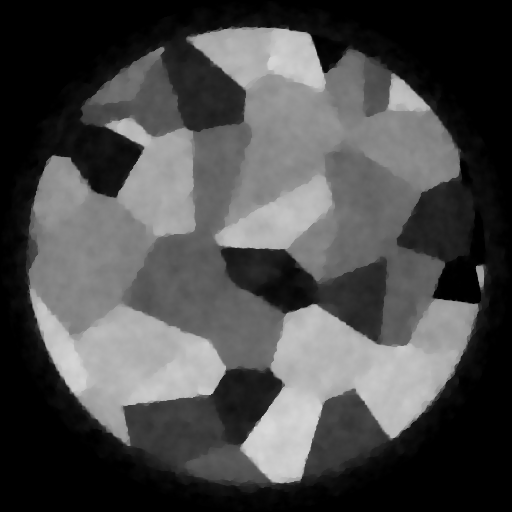}};
  \spy on (c5) in node at (4.3,-2.5);

  \node[label,inner sep=0pt] (c6) at (5,-2.2) {\includegraphics[width=0.32\columnwidth]{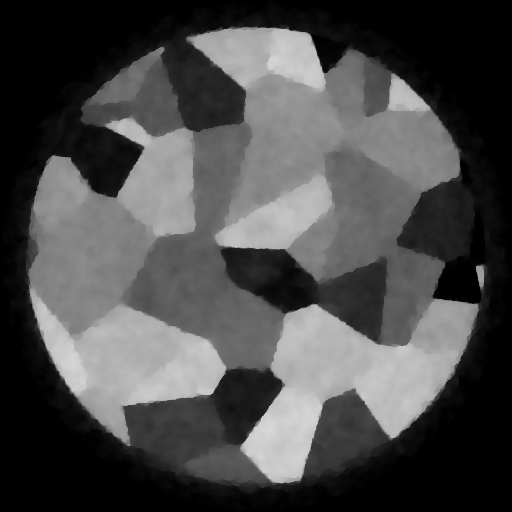}};
  \spy on (c6) in node at (5.3,-2.5);

\end{tikzpicture}
\caption{Phantom and reconstructions based on simulated low-intensity
  measurements. The display range for the images is 0 to 1.2
  cm$^{-1}$. The reconstructions with the TV-prior were computed with
  $\gamma = 3$. The insets are blow-ups of the reconstructions at the
  isocenter. \add{ The number of iterations was 500 for
    reconstructions without TV prior and 1,500 for reconstructions
  with TV prior.}}
\label{recnGP}
\end{figure*}

\newcommand{\hilite}{\cellcolor{blue!20}}
\begin{table*}
  \centering
\begin{tabular}{|l|d{3}{1}d{1}{2}d{1}{1}d{1}{2}|d{2}{1}d{1}{2}d{1}{1}d{1}{2}|d{2}{1}d{1}{2}d{1}{1}d{1}{2}|}
\hline
\multirow{2}{*}{Model}
  & \multicolumn{4}{c|}{Without TV (full domain, 2$\times$2 cm)}
  & \multicolumn{4}{c|}{Without TV (disc, radius 0.8 cm)}
  & \multicolumn{4}{c|}{TV $\gamma = 3$ (disc, radius 0.8 cm)} \\
	 & \multicolumn{1}{c}{RAE} &  \multicolumn{1}{c}{SSIM}
         & \multicolumn{1}{c}{RFE} &  \multicolumn{1}{c|}{RR}
         & \multicolumn{1}{c}{RAE} &  \multicolumn{1}{c}{SSIM}
         & \multicolumn{1}{c}{RFE} &  \multicolumn{1}{c|}{RR}
         & \multicolumn{1}{c}{RAE} &  \multicolumn{1}{c}{SSIM}
         & \multicolumn{1}{c}{RFE} &  \multicolumn{1}{c|}{RR} \\
 \hline \rowcolor[gray]{.8}
Baseline FBP 	 	 & 71.9 	 & 0.62 	 & 0.2 	 &    0.03 & 65.7 	 & 0.77      & 0.2 	 & 0.03  &
 \multicolumn{1}{c}{-}  & \multicolumn{1}{c}{-} & \multicolumn{1}{c}{-}&\multicolumn{1}{c|}{-}\\
FBP 	 	 & 101.2 	 & 0.55 	 & 3.5 	 & 0.66 	 & 94.7 	 & 0.70 	 & 3.6 	 & 0.66 &
\multicolumn{1}{c}{-}    & \multicolumn{1}{c}{-}&\multicolumn{1}{c}{-}&\multicolumn{1}{c|}{-}\\
P-FBP 	 	 & 73.1 	 & 0.62 	 & \hilite 1.8 	 & \hilite 0.20 	 & 66.7 	 & 0.77 	 & 1.6 	 & 0.19 &                                                                  \multicolumn{1}{c}{-} & \multicolumn{1}{c}{-}&\multicolumn{1}{c}{-}&\multicolumn{1}{c|}{-}\\
\hline \rowcolor[gray]{.8}
Baseline MAP 	 	 & 58.7 	 & 0.79 	 & 0.3 	 & 0.04 	 & 58.4 	 & 0.80 	 & 0.2 	 & 0.04 & 6.1 	 & 0.93 	 & 0.3 	 & 0.06 \\
AMAP 	 	 & 77.3 	 & 0.72 	 & 2.8 	 & 0.50 	 & 76.9 	 & 0.74 	 & 2.9 	 & 0.52 & 15.2 	 & 0.71 	 & 1.7 	 & 0.19 \\
WLS 	 	 & 76.9 	 & 0.72 	 & 2.8 	 & 0.50 	 & 76.6 	 & 0.74 	 & 2.9 	 & 0.52  & 15.2 	 & 0.71 	 & 1.7 	 & 0.19 \\
JMAP ($\beta=0$) 	 	 & 63.8 	 & 0.72 	 & 5.4 	 & 0.25 	 & 58.1 	 & \hilite 0.80 	 & 2.7 	 & \hilite 0.12 & 8.2 	 &\hilite 0.92 	 & 0.9 	 & \hilite 0.09 \\
SWLS ($\beta=0$) 	 	 & 63.9 	 & 0.72 	 & 5.5 	 & 0.26 	 & \hilite 58.0 	 & \hilite 0.80 	 & 2.7 	 &\hilite 0.12 & 8.3 	 &\hilite 0.92 	 & 1.0 	 & 0.10 \\
JMAP ($\beta=10$)	 	 & \hilite 61.8 	 & 0.74 	 & 3.1 	 & \hilite 0.20 	 & 58.4 	 & \hilite 0.80 	 &1.5 	 & 0.15  &\hilite 7.6 	 & \hilite 0.92 	 & \hilite 0.7 	 & \hilite 0.09 	\\
SWLS ($\beta=10$)  	 	 & \hilite 61.8 	 & 0.74 	 & 3.2 	 & \hilite 0.20 	 & 58.3 	 & \hilite 0.80 	 & \hilite 1.4 	 & 0.15 & 7.7 	 & \hilite 0.92 	 & 0.8 	 &  \hilite 0.09 \\
JMAP ($\beta=50$) 	 	 & 62.0 	 & \hilite 0.75 	 & 2.3 	 & 0.30 	 & 60.3 	 & 0.79 	 & 2.0 	 & 0.31 & \hilite 7.6 	 & 0.91 	 & 1.2 	 & 0.17 	\\
SWLS ($\beta=50$) 	 	 & 61.9 	 & \hilite 0.75 	 & 2.3 	 & 0.30 	 & 60.2 	 & 0.79 	 & 2.0 	 & 0.31 & 7.7 	 &0.91 	 & 1.3 	 & 0.17 \\
\hline
\end{tabular}
\medskip
  \caption{Error measures for reconstructions based on simulated low-intensity measurements.}
  \label{tab-ex1}
\end{table*}

In our first experiment, we used a phantom based on the ``grains''
phantom from AIR Tools, shown in the upper left corner of Fig.~\ref{recnGP}. We applied a
circular mask of radius 0.8 cm to obtain a phantom that is fully
contained by the reconstruction grid. We used $I_0=500$ in this
experiment, corresponding to approximately 500 photons per detector
element per projection. As a result, the SNR is relatively
low. Estimates based on low SNR measurements generally have a high
variance, and hence a good model and strong priors are of paramount
importance. The reconstructions shown in Fig.~\ref{recnGP} demonstrate
this. The baseline reconstructions were computed using the true
flat-field, and hence they are ``inverse crime'' reconstructions that
serve only as a baseline for comparison. The two baseline MAP
reconstructions (with and without the TV prior) are based on the model
\eqref{e-map-map}. Using the flat-field estimate $\vml$ instead of the
true flat-field, we obtained the FBP and AMAP reconstructions. It is
clear from these reconstructions that the flat-field estimation errors
introduce severe ring artifacts, even in the presence of a strong
prior such as the TV-prior. The ring artifacts are especially severe
near the center of the image (cf.\ Section \ref{effecffferror}).

The preprocessed FBP (P-FBP) reconstruction is the
result of applying the combined wavelet and FFT filtering
preprocessing method\footnote{We used a damping factor of 0.9 and a
  Daubechies 5 wavelet with a three-level decomposition.} by M{\"u}nch
et al.\ \cite{Munch2009} to the sinogram, followed by FBP. This removes stripes from
the sinogram, and although there are still some noticeable ring
artifacts in the reconstruction, the preprocessing step clearly reduces the severity of
the artifacts. However, the preprocessing step involves several
parameters that must be carefully tuned, and it does not directly
allow us to use the AMAP or MAP-based reconstruction models for
reconstruction.

The proposed models are quite effective at reducing ring artifacts, as
can be seen from the JMAP reconstructions as well as the SWLS
reconstruction. Notice that both the SWLS ($\beta=0$) reconstruction and the  JMAP
($\beta=0$) reconstruction without the TV prior do not involve any parameters.

\begin{figure*} 
  \centering
  \begin{tikzpicture}[font=\footnotesize]
    \begin{semilogyaxis}[name=plot1,yscale=0.45, no markers,xmin=0,xmax=1400,ymin=3e0,ymax=2e2,width=0.95\columnwidth,
      ylabel={Relative attenuation error},xlabel={},grid=major,legend
      columns=2,legend style={font=\scriptsize,
        draw=none,anchor=west,column sep=4mm},
      legend entries={Baseline MAP,Baseline MAP (init.~FBP),AMAP,AMAP
        (init.~P-FBP),JMAP,JMAP (init.~P-FBP),JMAP (init.~AMAP)},
      transpose legend,
      legend cell align=left,
      legend to name=named1,
      title={Without TV-regularization},
      title style={at={(0.5,2.3)},anchor=south}]

      \addplot[thick, red!65] table [x index=0, y index=1]{data/ex2_RAE_notv.dat};
      \addplot[thick, red!65, dashed] table [x index=0, y index=2]{data/ex2_RAE_notv.dat};

      \addplot[thick, blue!35] table [x index=0, y index=3]{data/ex2_RAE_notv.dat};
      \addplot[thick, blue!35, dashed] table [x index=0, y index=4]{data/ex2_RAE_notv.dat};

      \addplot[thick, black] table [x index=0, y index=5]{data/ex2_RAE_notv.dat};
      \addplot[thick, black, dashed] table [x index=0, y index=6]{data/ex2_RAE_notv.dat};
      \addplot[thick, black, dash dot] table [x index=0, y index=7]{data/ex2_RAE_notv.dat};

    \end{semilogyaxis}
    \begin{semilogyaxis}[name=plot2, at = (plot1.right of east), anchor= right of east,yscale=0.45, no markers,xmin=0,xmax=1400,ymin=3e0,ymax=2e2,width=0.95\columnwidth,
      ylabel={Relative attenuation error},xlabel={},grid=major,xshift =\columnwidth,
      title={With TV-regularization ($\gamma =3$)},
      title style={at={(0.5,2.3)},anchor=south}
      ]
      \addplot[thick, red!65] table [x index=0, y index=1]{data/ex2_RAE_tv.dat};
      \addplot[thick, red!65, dashed] table [x index=0, y index=2]{data/ex2_RAE_tv.dat};

      \addplot[thick, blue!35] table [x index=0, y index=3]{data/ex2_RAE_tv.dat};
      \addplot[thick, blue!35, dashed] table [x index=0, y index=4]{data/ex2_RAE_tv.dat};

      \addplot[thick, black] table [x index=0, y index=5]{data/ex2_RAE_tv.dat};
      \addplot[thick, black, dashed] table [x index=0, y index=6]{data/ex2_RAE_tv.dat};
      \addplot[thick, black, dash dot] table [x index=0, y index=7]{data/ex2_RAE_tv.dat};

    \end{semilogyaxis}
     \begin{semilogyaxis}[name=plot3, at = (plot1.below south west), anchor=north west,yscale=0.45, no markers,xmin=0,xmax=1400,ymin=1e-2,ymax=1e1,width=0.95\columnwidth,
      ylabel={Ring ratio},xlabel={Iterations},grid=major
      ]
      \addplot[thick, red!65] table [x index=0, y index=1]{data/ex2_RR_notv.dat};
      \addplot[thick, red!65, dashed] table [x index=0, y index=2]{data/ex2_RR_notv.dat};

      \addplot[thick, blue!35] table [x index=0, y index=3]{data/ex2_RR_notv.dat};
      \addplot[thick, blue!35, dashed] table [x index=0, y index=4]{data/ex2_RR_notv.dat};

      \addplot[thick, black] table [x index=0, y index=5]{data/ex2_RR_notv.dat};
      \addplot[thick, black, dashed] table [x index=0, y index=6]{data/ex2_RR_notv.dat};
      \addplot[thick, black, dash dot] table [x index=0, y index=7]{data/ex2_RR_notv.dat};

      \addplot[color=darkgray,loosely dashed] coordinates {(50,0.01)(50,10)};
    \end{semilogyaxis}
    \begin{semilogyaxis}[name=plot4, at = (plot2.below south west), anchor=north west ,yscale=0.45, no markers,xmin=0,xmax=1400,ymin=1e-2,ymax=1e1,width=0.95\columnwidth,
      ylabel={Ring ratio},xlabel={Iterations},grid=major
      ]
      \addplot[thick, red!65] table [x index=0, y index=1]{data/ex2_RR_tv.dat};
      \addplot[thick, red!65, dashed] table [x index=0, y index=2]{data/ex2_RR_tv.dat};

      \addplot[thick, blue!35] table [x index=0, y index=3]{data/ex2_RR_tv.dat};
      \addplot[thick, blue!35, dashed] table [x index=0, y index=4]{data/ex2_RR_tv.dat};

      \addplot[thick, black, ] table [x index=0, y index=5]{data/ex2_RR_tv.dat};
      \addplot[thick, black, dashed] table [x index=0, y index=6]{data/ex2_RR_tv.dat};
      \addplot[thick, black, dash dot] table [x index=0, y index=7]{data/ex2_RR_tv.dat};

       \addplot[color=darkgray,loosely dashed] coordinates {(50,0.01) (50,10)};
    \end{semilogyaxis}
\end{tikzpicture}
\\
\ref{named1}
\caption{Results of semi-convergence and initialization study. \add{Reconstructions are computed with a UP prior $\beta = 0$.}}
\label{ex2}
\end{figure*}
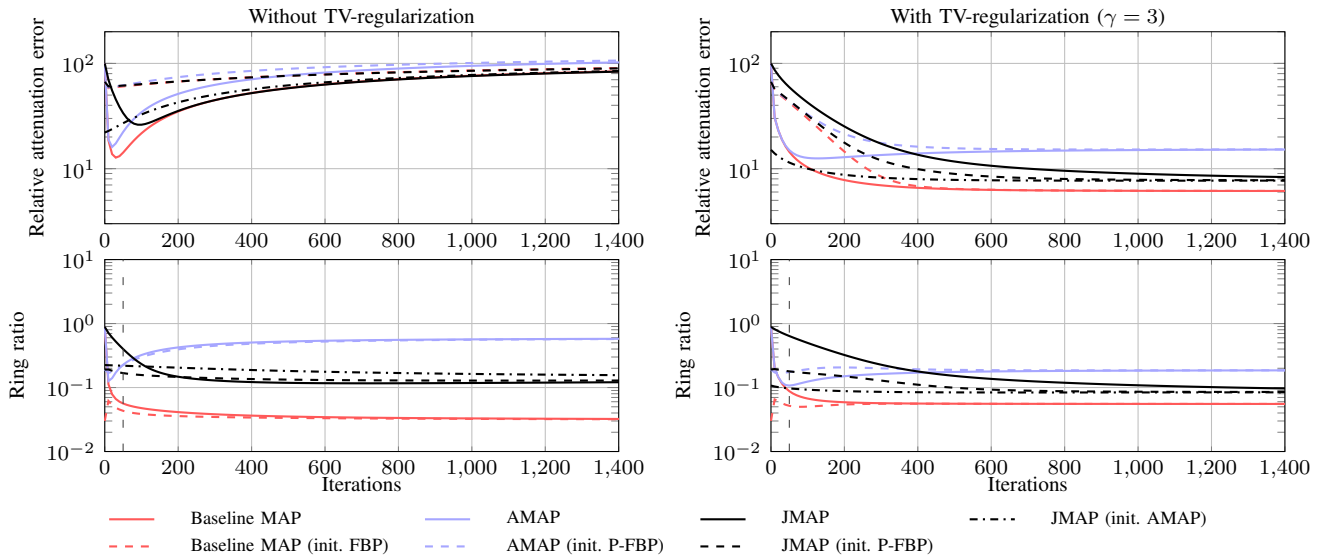

For the experiments without the TV-prior, Table \ref{tab-ex1} shows the
error measures based on both the full reconstruction domain and based
on a disc of radius 0.8 cm (corresponding to the support of the
phantom). The latter approach ignores noise and ring artifacts outside
the phantom, and hence this gives a more practical picture of the
performance. For the reconstructions with a TV-prior, we report our
results based on a disc of radius 0.8 cm. Notice that in all cases, we
obtain the best reconstruction (in terms of both RAE and SSIM) using
either the JMAP reconstruction model or the SWLS model. Moreover,
these reconstructions have RAEs that are similar to those of the
baseline MAP reconstructions. We also see that RRs and the RAEs for
the JMAP reconstructions appear to be correlated, but interestingly,
the RFEs do not seem to agree with the RAEs.

Despite the fact that the P-FBP reconstruction is worse than the JMAP
reconstructions, it is interesting to note that it may be used to
compute an improved flat-field estimate. In our experiment, the ML
estimate $\vml$ had a relative error of 4.8\%, but the flat-field
estimate computed based on the P-FBP reconstruction had a relative
error of only around 1.8\%. However, using the TV-prior, the JMAP
and SWLS model still produced the best flat-field estimate of all the models.

Finally, we remark that the AMAP and WLS reconstructions may be
improved slightly by increasing the parameter $\gamma$. Using
$\gamma = 10$, we obtained AMAP and WLS reconstructions with a
relative error of around 10\%, and although these reconstructions did
not have noticeable ring artifacts, they contained an increased amount
of undesirable TV-artifacts. On the other hand, the JMAP and SWLS
reconstructions obtained with $\gamma = 3$ only have a limited amount
of ring artifacts and TV artifacts, and hence we conclude that the
proposed model allows us to reduce ring artifacts using a smaller
regularization parameter $\gamma$ than with the AMAP or WLS models,
thus limiting unnecessary TV-induced artifacts.

\subsubsection{Semi-convergence and Initialization}

We now investigate the role of regularization and its influence on the
reconstruction. Recall that X-ray tomographic imaging is an ill-posed
problem where a small amount of noise in the measurements may results
in a large change in the reconstruction if it is not regularized by a
suitable prior. Thus, without regularization, intermediate iterates
sometimes provide better reconstructions than iterates close to
convergence. This behavior is known as semi-convergence and depends
on the reconstruction method as well as initialization. Semi-convergence
behavior often indicates that the reconstruction is under-regularized, and
hence a solution to our convex reconstruction model may be a poor
reconstruction. In practice it is difficult to rely on
semi-convergence as the true solution is unknown.

We use the same experimental setup as in the previous experiment\remove{, with
the exception that $I_0 = 5,\!000$, corresponding to a small increase in
intensity}.  Fig.~\ref{ex2} shows RAE and RR as a function of the
number of iterations, with and without the TV-prior (\ie,
regularization). The semi-convergence behavior is evident without the
TV-prior, and not surprisingly, the baseline reconstruction obtains the lowest RAE at the
semi-convergence point after approximately \remove{$100$} \add{$50$} iterations. After the
semi-convergence point, noise start to dominate the reconstruction
and the RAE starts to increase monotonically. Comparing the AMAP and
JMAP models, we see that the AMAP model has a lower RAE at the
semi-convergence point, but it converges to a higher
RAE. Taking the definition of the AMAP and JMAP estimators into
account, we can conclude that the JMAP model still converges to a
better reconstruction than the AMAP model. Fig.~\ref{ex2} also shows
the RR error measure, and while the AMAP model exhibits
semi-convergence both with respect to the RAE and the RR, the JMAP model
appears to monotonically reduce the RR despite semi-convergence with
respect to the RAE.

The dashed curves in Fig.~\ref{ex2} show the results of the same experiment, but using
the P-FBP reconstruction of $u$ as initialization (the baseline MAP
was initialized with the baseline FBP reconstruction). The FBP
reconstruction has a smaller RAE than the zero-initialization, but FBP
reconstructions may be quite noisy when the SNR is low. Consequently,
this initialization may not lead to faster convergence without
regularization, as can be seen in Fig \ref{ex2}. The figure also shows
that the AMAP reconstruction method still exhibits a mild degree of
semi-convergence when using the TV-prior, but the baseline method and
the JMAP method appear to reduce the RAE and the RR
monotonically. Moreover, it is clear that the FBP-initialization helps
when combined with the TV-prior. Finally, using the 50th AMAP iterate
as initialization for JMAP (corresponding to the semi-convergence
point for the RR), we obtained a significant improvement in the number
of iterations when compared to initialization with zeros.

\subsubsection{Noise Analysis}
\begin{figure*}
  \centering
  \begin{tikzpicture}[x=0.39\columnwidth,y=0.39\columnwidth,font=\footnotesize,spy scope={magnification=2, size=0.18\columnwidth,every spy in node/.style={draw}}]
    \node[label={above:Baseline},label={[rotate=90,anchor=south]left:Bias}]
    (a1) at (0,0) {\includegraphics[width=0.38\columnwidth]{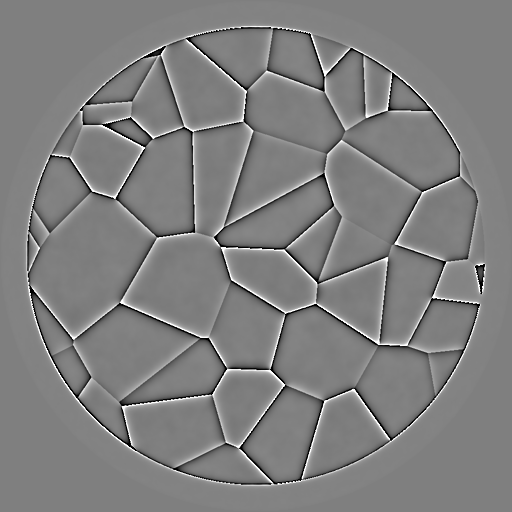}};
  \node[label={above:AMAP}] (a2) at (1,0) {\includegraphics[width=0.38\columnwidth]{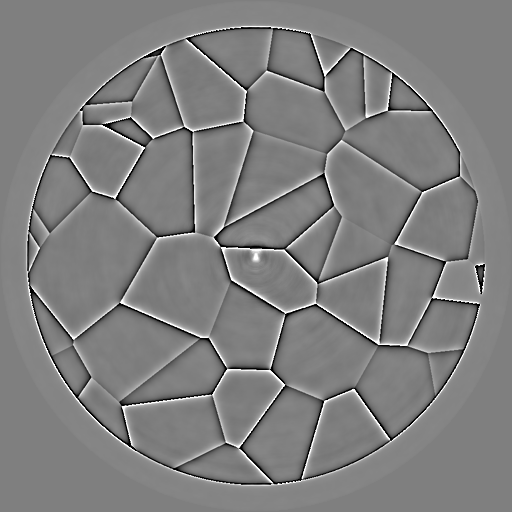}};
  \node[label={above:JMAP ($\beta=0$)}] (a3) at (2,0)
  {\includegraphics[width=0.38\columnwidth]{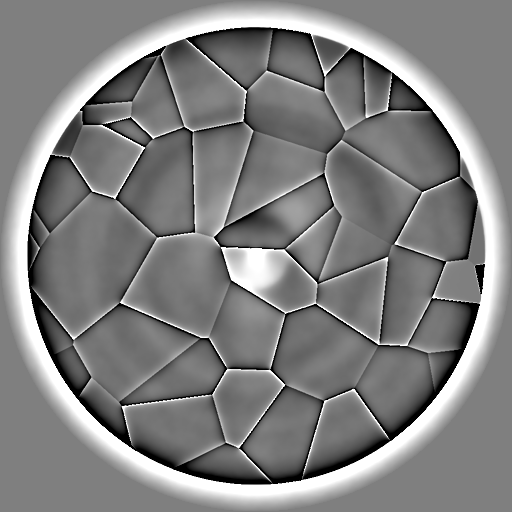}};
  \node[label={above:JMAP ($\beta=10$)}] (a4) at (3,0)
  {\includegraphics[width=0.38\columnwidth]{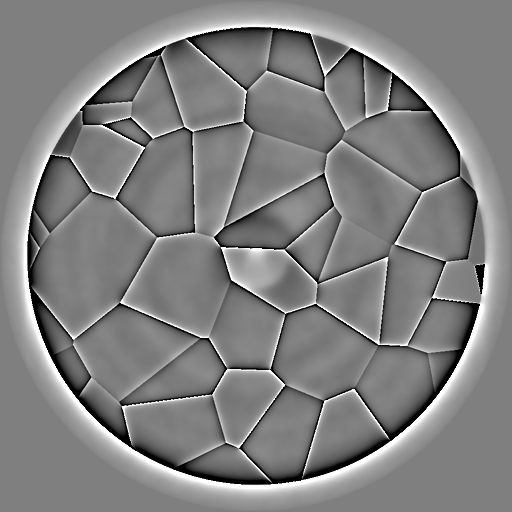}};
  \node[label={above:JMAP ($\beta=50$)}] (a5) at (4,0)
  {\includegraphics[width=0.38\columnwidth]{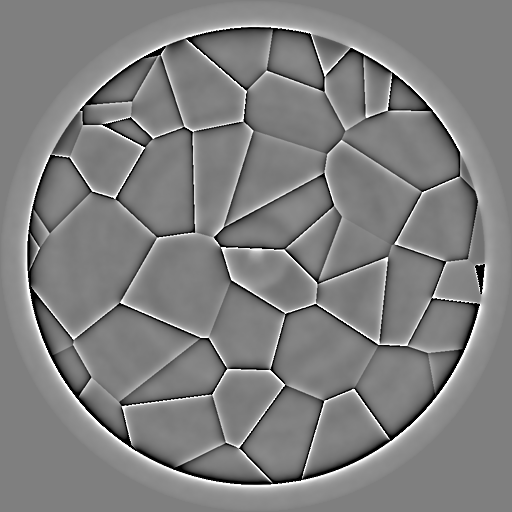}};
  \spy on (a1) in node at (-0.25,0.25);
  \spy on (a2) in node at (0.75,0.25);
  \spy on (a3) in node at (1.75,0.25);
  \spy on (a4) in node at (2.75,0.25);
  \spy on (a5) in node at (3.75,0.25);

  \node[label={[rotate=90,anchor=south]left:Std.~deviation}] (b1) at (0,-1) {\includegraphics[width=0.38\columnwidth]{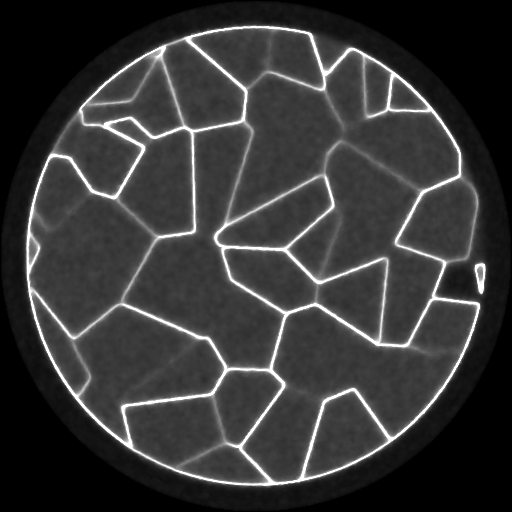}};
  \node (b2) at (1,-1) {\includegraphics[width=0.38\columnwidth]{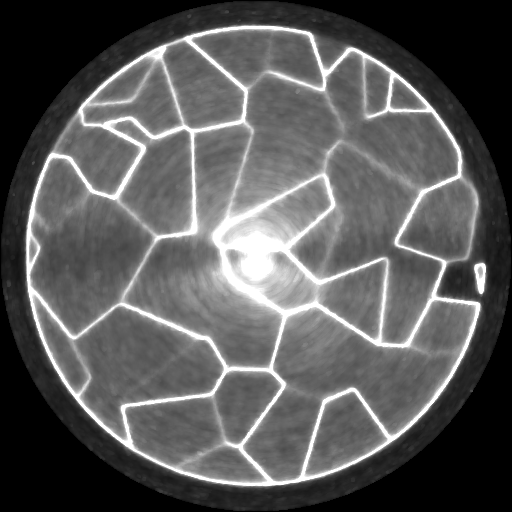}};
  \node (b3) at (2,-1) {\includegraphics[width=0.38\columnwidth]{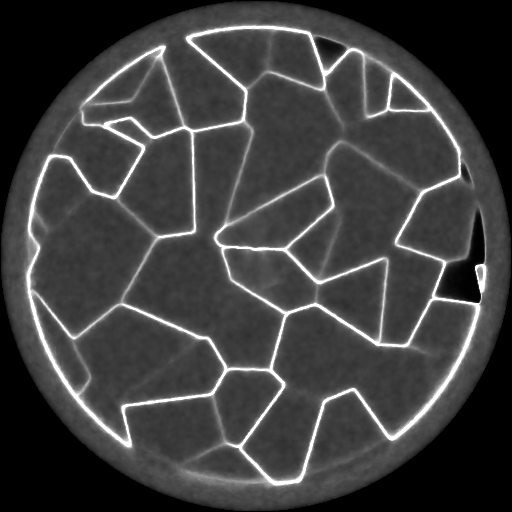}};
  \node (b4) at (3,-1) {\includegraphics[width=0.38\columnwidth]{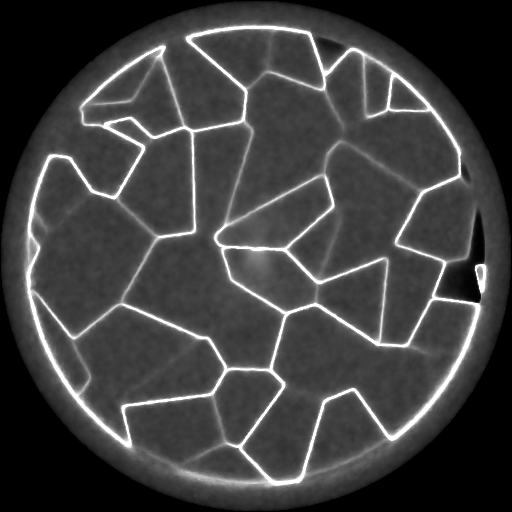}};
  \node (b5) at (4,-1) {\includegraphics[width=0.38\columnwidth]{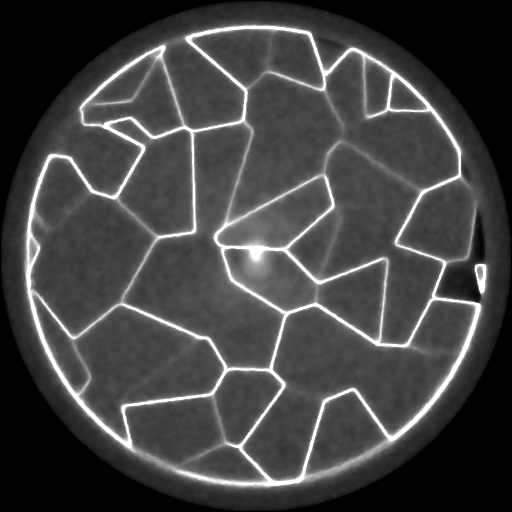}};
  \spy on (b1) in node at (-0.25,-0.75);
  \spy on (b2) in node at (0.75,-0.75);
  \spy on (b3) in node at (1.75,-0.75);
  \spy on (b4) in node at (2.75,-0.75);
  \spy on (b5) in node at (3.75,-0.75);
  \end{tikzpicture}
  \caption{Pixelwise bias and standard deviation based on 200
    realizations of all measurements. The
    display range for the bias images is $-$0.1 to 0.1 cm$^{-1}$, and
    the display range for the standard deviation images is 0 to 0.06 cm$^{-1}$.  The reconstructions are computed with TV-prior with $\gamma=3$. The insets are blow-ups of the reconstructions at the isocenter.}
  \label{fig-noiseanalysis}
\end{figure*}

To investigate the noise properties of the proposed reconstruction
model, we generated 200 realizations of all measurements based on the
grains phantom (see Fig.~\ref{recnGP}) and with $I_0=500$. We then
computed pixelwise bias (the difference between the mean of the
reconstructions and the phantom) and standard deviation for
reconstructions based on the baseline MAP, the AMAP, and the JMAP
reconstruction models. All reconstructions were computed with the
TV-prior ($\gamma=3$) and 1,500 iterations. The results are shown in
Fig.~\ref{fig-noiseanalysis}. Generally speaking, the AMAP model is
less biased than the JMAP model. For small values of $\beta$, the JMAP
bias is somewhat large in comparison to the AMAP bias, especially near
the boundary of the object and at the isocenter. However, the JMAP
bias decreases when the parameter $\beta$ is increased, but at the
cost of increasing the standard deviation. This is consistent with the
fundamental trade-off between bias and variance in statistical
learning.  More importantly, the standard deviation is significantly
lower for the JMAP model in comparison to the AMAP model, and it is
even comparable to that of the baseline MAP model when $\beta$ is
small. Notice that in all instances, the standard deviation is
particularly large near the interfaces of the grains where the
intensity jumps.

Recall from the previous experiment that the flat-field estimate may
converge very slowly. As a consequence, the bias component that is
induced by flat-field estimation errors decreases slowly as we
increase the number of iterations. The results therefore depend on the
stopping criteria (\ie, the number of iterations).  Finally we note
that the noise results for the SWLS model were very similar to those
of the the JMAP model, and hence we have chosen to omit the SWLS
results for the sake of brevity.

\subsubsection{Flat-field Regularization}

\begin{figure}
\centering

\begin{tikzpicture}[x=0.305\columnwidth,y=0.305\columnwidth,font=\footnotesize]
  \node[label={above:AMAP},label={[rotate=90,anchor=south]left:Reconstruction}] at (0,0) {\includegraphics[width=0.3\columnwidth]{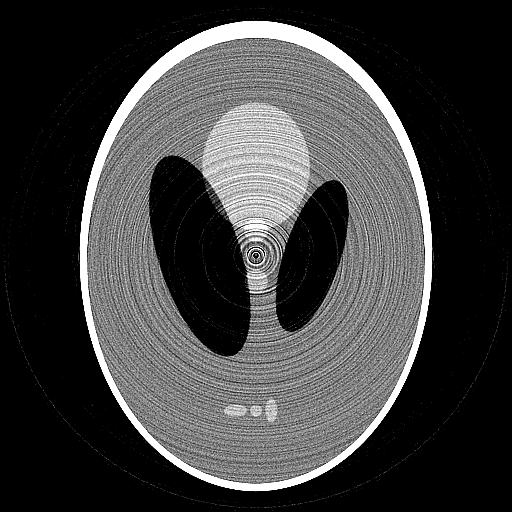}};
  \node[label={above:JMAP ($\beta=0$)}] at (1,0) {\includegraphics[width=0.3\columnwidth]{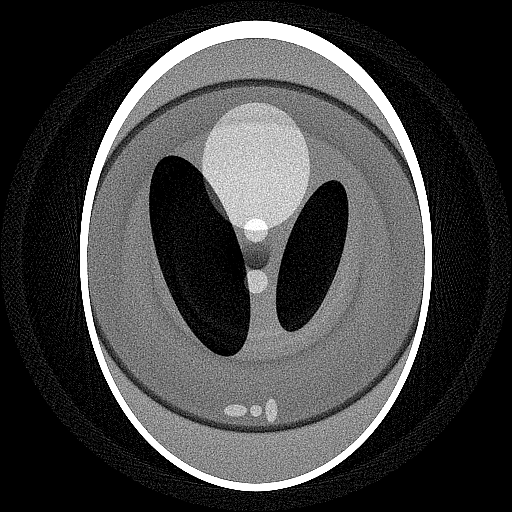}};
  \node[label={above:JMAP ($\beta=50$)}] at (2,0) {\includegraphics[width=0.3\columnwidth]{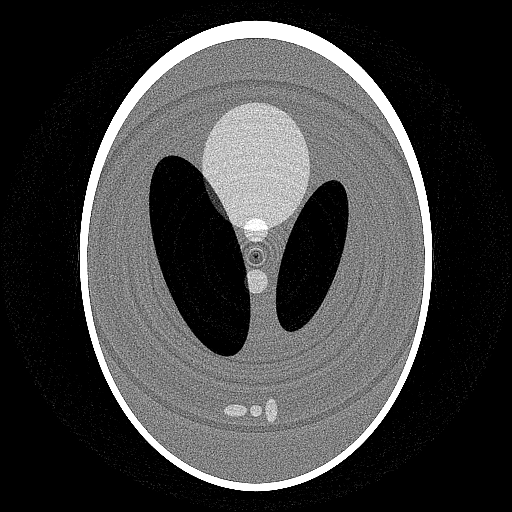}};

  \node[label={[rotate=90,anchor=south]left:$|\psi_v(\hat v)|$}] at (0,-1) {\includegraphics[width=0.3\columnwidth]{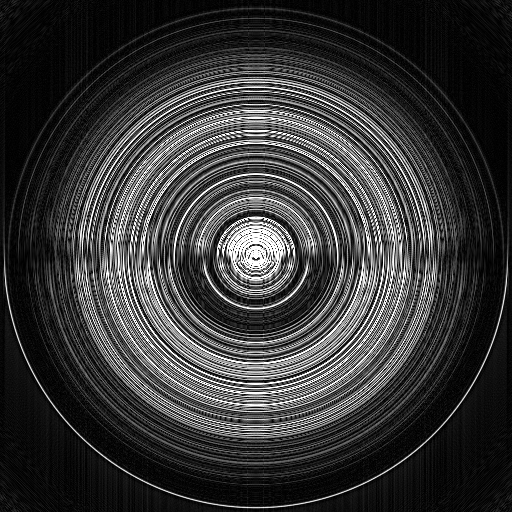}};
  \node at (1,-1) {\includegraphics[width=0.3\columnwidth]{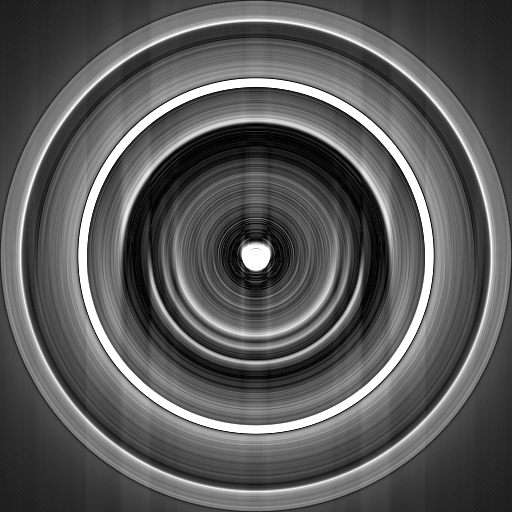}};
  \node at (2,-1) {\includegraphics[width=0.3\columnwidth]{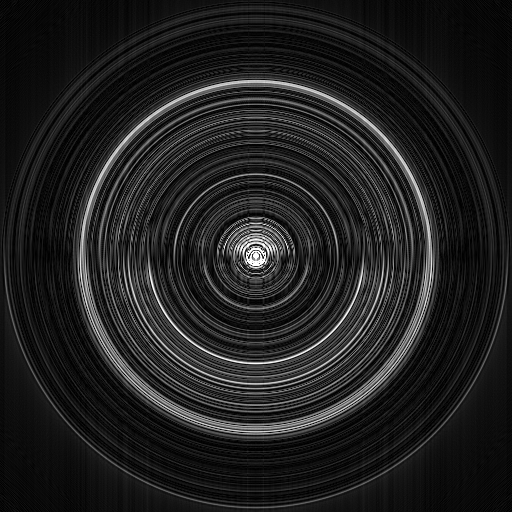}};
\end{tikzpicture}

\par\medskip
\begin{tikzpicture}[font=\footnotesize]
    \begin{axis}[name=plot1,title={$\beta=0$},title style={yshift=-0.5em},y post scale=1.0,no markers,xmin=1,xmax=512,ymin=0e0,ymax=1e0,width=0.58\columnwidth,
      ylabel={},xlabel={Detector element},label style={yshift=0.4em},legend
      columns=-1,legend style={draw=none, column sep=2mm,font=\footnotesize},
      legend entries={$\theta_1 (\vml)$,$\theta_2 (\vmly)$,$\theta_3 (\hat{v}_{\mathrm{pr}})$},
      legend to name=named2,xtick={1,256,512}]
      \addplot[thick, blue, dashed] table [x index=0, y index=1]{data/ex3_beta_jmap.dat};
      \addplot[thick, black, dashdotted] table [x index=0, y index=2]{data/ex3_beta_jmap.dat};
      \addplot[thick, red] table [x index=0, y index=3]{data/ex3_beta_jmap.dat};
\end{axis}
    \begin{axis}[name=plot2, title={$\beta=50$},title style={yshift=-0.5em},at = (plot1.right of east), anchor= right of east,xshift=0.5\columnwidth,y post scale=1.0, xscale=1,no markers,xmin=1,xmax=512,ymin=0e0,ymax=1e0,width=0.58\columnwidth,
      ylabel={},xlabel={Detector element},label style={yshift=0.4em},xtick={1,256,512}]
      \addplot[thick, blue, dashed] table [x index=0, y index=1]{data/ex3_beta_jmap2.dat};
      \addplot[thick, black, dashdotted] table [x index=0, y index=2]{data/ex3_beta_jmap2.dat};
      \addplot[thick, red] table [x index=0, y index=3]{data/ex3_beta_jmap2.dat};
\end{axis}
\end{tikzpicture}
\ref{named2}
\par\medskip
\begin{tikzpicture}[font=\footnotesize]
    \begin{axis}[name=plot3,title={},title style={yshift=-0.5em},y post scale=0.4,no markers,xmin=1,xmax=512,width=\columnwidth,
      ylabel={Rel.~error (\%)},xlabel={Detector element},label style={yshift=-0.4em},legend
      columns=0,legend style={draw=none, column sep=2mm,font=\footnotesize},
      legend entries={$e_{v}(\vml)$,$e_{v}(\hat v(u))\, (\beta = 0)$,$e_{v}(\hat v(u))\, (\beta = 50)$},
      legend to name=named3,xtick={1,256,512}]
      \addplot[thin, blue, densely dotted] table [x index=0, y index=1]{data/ex3_vhu.dat};
      \addplot[thick, black, densely dashed ] table [x index=0, y index=2]{data/ex3_vhu.dat};
      \addplot[thick, red] table [x index=0, y index=3]{data/ex3_vhu.dat};
\end{axis}
\end{tikzpicture}
\ref{named3}
\par\medskip
{\footnotesize
\begin{tabular}{|l|d{3}{1}d{1}{2}d{1}{1}d{1}{2}|}
\hline
\multirow{1}{*}{Model}
	 & \multicolumn{1}{c}{Baseline} &  \multicolumn{1}{c}{AMAP}
         & \multicolumn{1}{c}{JMAP ($\beta = 0$)} &  \multicolumn{1}{c|}{JMAP ($\beta = 50$)}\\
 \hline
RAE	 	 & 13.49 	 & 17.37 	 & 17.56	 &    14.00 \\
RR	 	 & 0.03 	 & 0.60 	 & 0.96 	 & 0.31 \\
\hline
\end{tabular}
}
\par\medskip
\caption{Reconstructions of the Shepp--Logan phantom after 1,000
  iterations, without the  TV-prior on $u$. The display range for the
  reconstruction images is 0 to 0.4 cm$^{-1}$, and 0 to 0.04 for the
  ring images $\psi_v({\hat{v}})$. \remove{The two plots show the values of
  $\theta_1,\theta_2,\theta_3$, as defined in \eqref{theta-par}, for
  $\beta=0$ and $\beta=50$, and the table lists the RAE and RR error
  measures.}\add{The first two plots show the values of
  $\theta_1,\theta_2,\theta_3$, as defined in \eqref{theta-par}, for
  $\beta=0$ and $\beta=50$. The third plot shows the element-wise
  relative error with respect to true flat-field $v$, defined as $e_{v}(\hat v) = 100\cdot\diag(v)^{-1}(\hat v - v)$, for the ML flat-field
  estimate $\vml$ and two JMAP flat-field estimates.
 The table lists the RAE and RR error
  measures.}}
\label{fig-recnSL}
\end{figure}

Our next experiment demonstrates a potential shortcoming of the
proposed model when using the UP flat-field prior for
reconstruction. We used the Shepp--Logan phantom for the experiment,
but unlike in the previous experiments, we generated the measurements
by evaluating the line integrals analytically. The intensity parameter
was $I_0= 10^5$. The reconstruction based on \eqref{e-map-amap}, the
leftmost reconstruction in Fig.~\ref{fig-recnSL}, has some low-level
ring artifacts. JMAP with the FE prior and $\beta=0$ leads to the
reconstruction in the middle of Fig.~\ref{fig-recnSL}. Somewhat
surprisingly, while the low-level rings are mostly gone, the
reconstruction has a few wide and very noticeable rings. These rings
arise because of the structure of the flat-field estimation errors
which can be seen by looking at the reconstruction $\psi_v(\hat v)$,
defined in \eqref{e-ring-err} and shown in
Fig.~\ref{fig-recnSL}. Several high-intensity rings appear clearly,
and these can be linked to large flat-field estimation errors
associated with a small number of detector elements. In particular,
the detector elements corresponding to rays that intersect the outer
ellipsoidal shell of the Shepp--Logan phantom tangentially give rise
to large estimation errors. \add{We remark that we have observed
  experimentally that these artifacts seem to be exacerbated by the
  fact that the two outer Shepp--Logan ellipses are centered at the isocenter.}

Now recall that the flat-field estimate $\hat{v}(u)$ can be expressed
as \eqref{theta-par}, \ie, a convex combination of independent
estimates. Thus, the weights $\theta$ indicate the emphasis of the
different flat-field estimates. The plots in Fig.~\ref{fig-recnSL}
show these weights for two different priors parameterized by
$\beta$. We see that when $\beta=0$ (corresponding to the UP
flat-field prior), the flat-field estimate is based
almost entirely on $\vmly$, and the estimates $\vml$ and
$\hat v_{\mathrm{pr}}$ both receive \remove{a negligible
  weight }\add{negligible (but nonzero) weights}. \add{Inspecting the
  corresponding flat-field estimate (the bottom plot in
  Fig.~\ref{fig-recnSL}) reveals that for $\beta = 0$, the JMAP estimate is worse
  than the ML estimate $\vml$. This indicates over-fitting.} To
mitigate this, we can emphasize the flat-field ML estimate $\vml$ by
using the FE prior (\ie, $\alpha = 1+ \beta(\vml)$), as described in
\ref{ss-hyperparameters}. Doing so effectively removes the major rings
that were present with the \remove{UP prior
($\beta=0$)}\add{FE prior with $\beta = 0$}, as shown in the rightmost reconstruction in
Fig.~\ref{fig-recnSL}. Moreover, the rightmost plot in the figure
confirms that the resulting
flat-field estimate depends less on $\vmly$ than with the \remove{UP
  prior}\add{FE prior with $\beta =0$}. The
FBP reconstructions of the flat-field error, shown below the
reconstructions in Fig. \ref{fig-recnSL}, clearly show a reduction in
ring artifacts compared to the basic AMAP and JMAP
reconstructions. 

\subsection{Real Data Study}

\begin{figure*}
\centering
 \begin{tikzpicture}[x=0.20\textwidth,y=0.20\textwidth,font=\footnotesize,spy scope={magnification=\magn, size=0.0875\textwidth,every spy in node/.style={draw}}]
  \node[label={[font=\footnotesize]below:FBP},inner sep=0pt] (a1) at (0,0) {\includegraphics[width=0.19\textwidth]{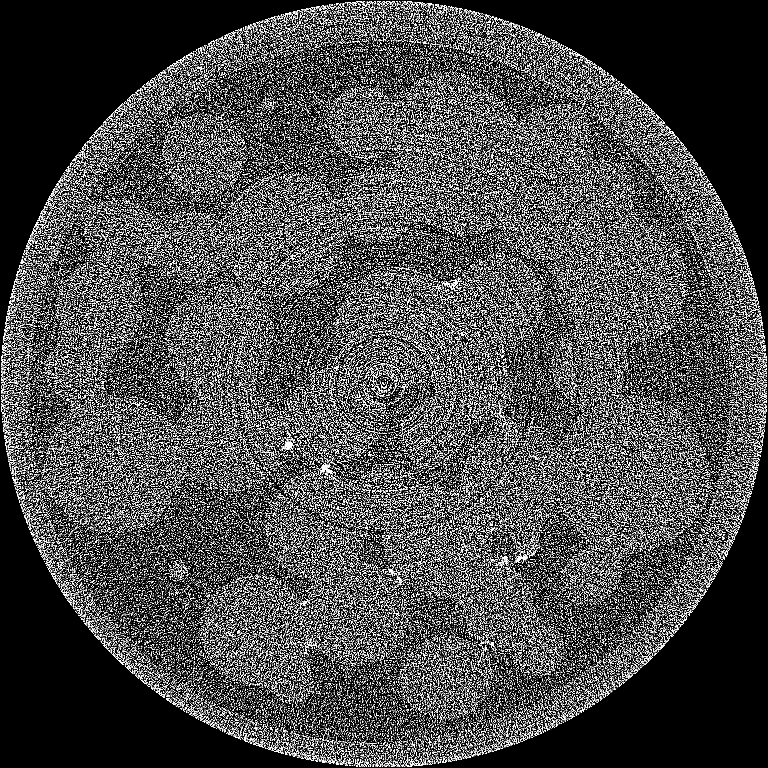}};
  \spy on (a1) in node at (-0.255,0.255);
  \node[label={[font=\footnotesize]below:P-FBP},inner sep=0pt] (a2) at (1,0) {\includegraphics[width=0.19\textwidth]{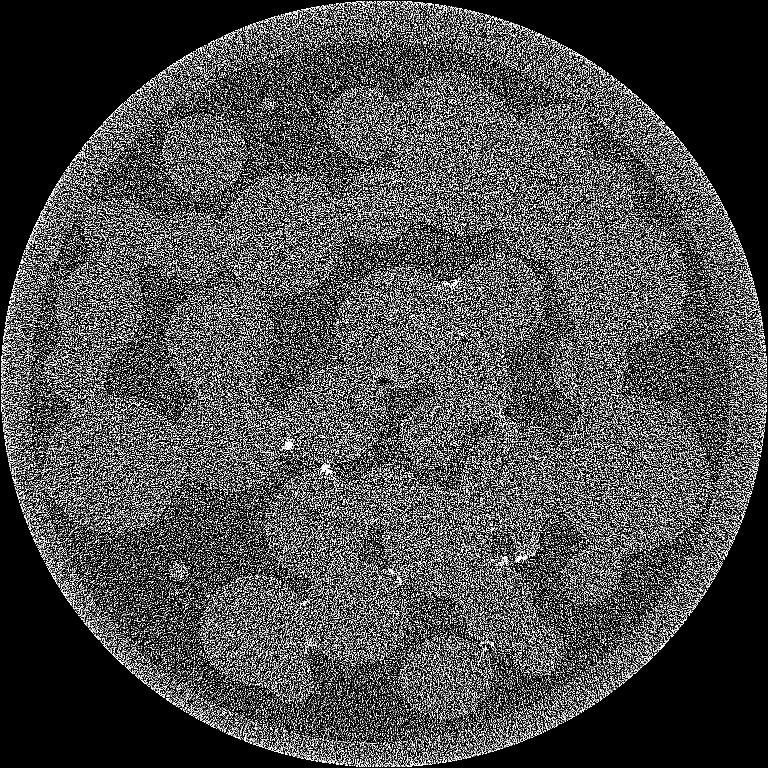}};
  \spy on (a2) in node at (0.745,0.255);
  \node[label={[font=\footnotesize]below:AMAP},inner sep=0pt] (a3) at (2,0) {\includegraphics[width=0.19\textwidth]{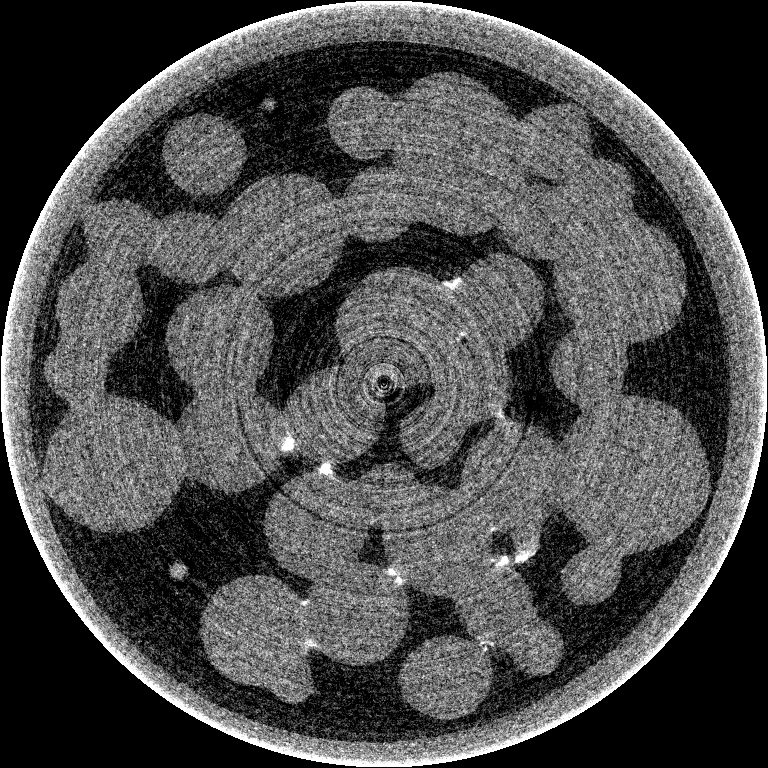}};
  \spy on (a3) in node at (1.745,0.255);
  \node[label={[font=\footnotesize]below:JMAP ($\beta=0$)},inner sep=0pt] (a4) at (3,0) {\includegraphics[width=0.19\textwidth]{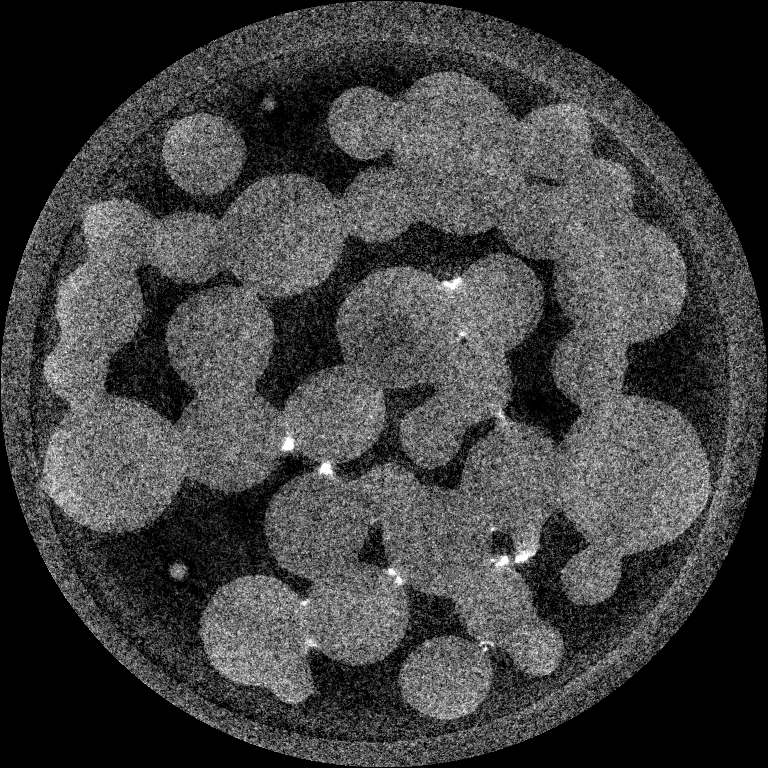}};
  \spy on (a4) in node at (2.745,0.255);
  \node[label={[font=\footnotesize]below:JMAP ($\beta=200$)},inner sep=0pt] (a5) at (4,0) {\includegraphics[width=0.19\textwidth]{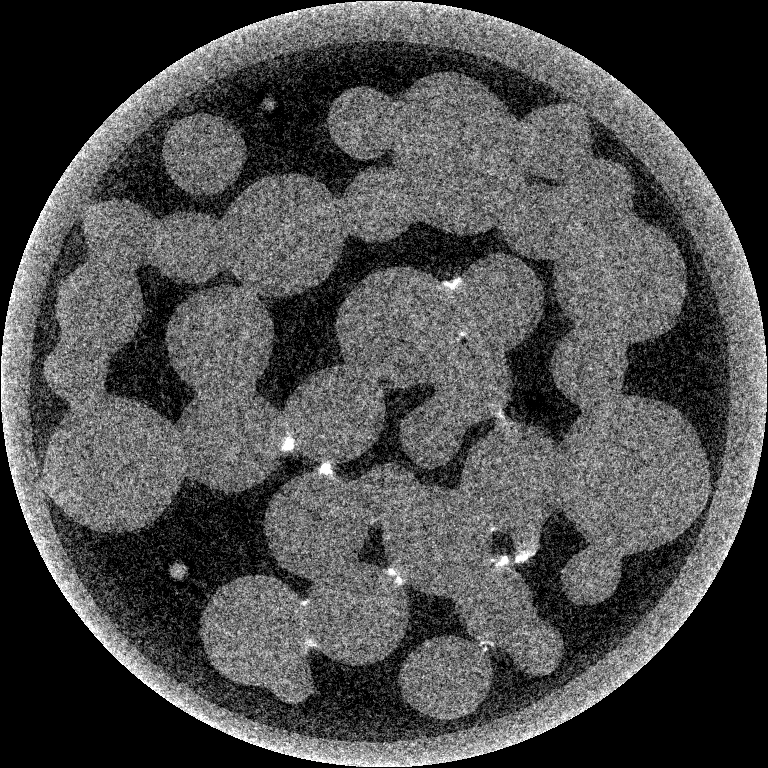}};
  \spy on (a5) in node at (3.745,0.255);

  \node[label={[font=\footnotesize]below:FBP + smoothing},inner sep=0pt] (b1) at (0,-1.1) {\includegraphics[width=0.19\textwidth]{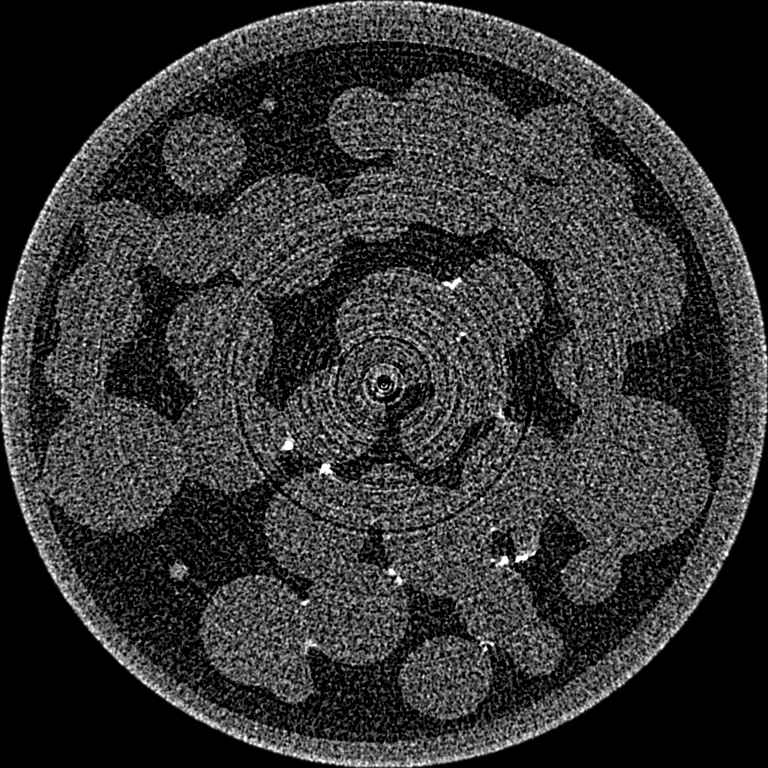}};
  \spy on (b1) in node at (-0.255,-0.845);
  \node[label={[font=\footnotesize]below:P-FBP + smoothing},inner sep=0pt] (b2) at (1,-1.1) {\includegraphics[width=0.19\textwidth]{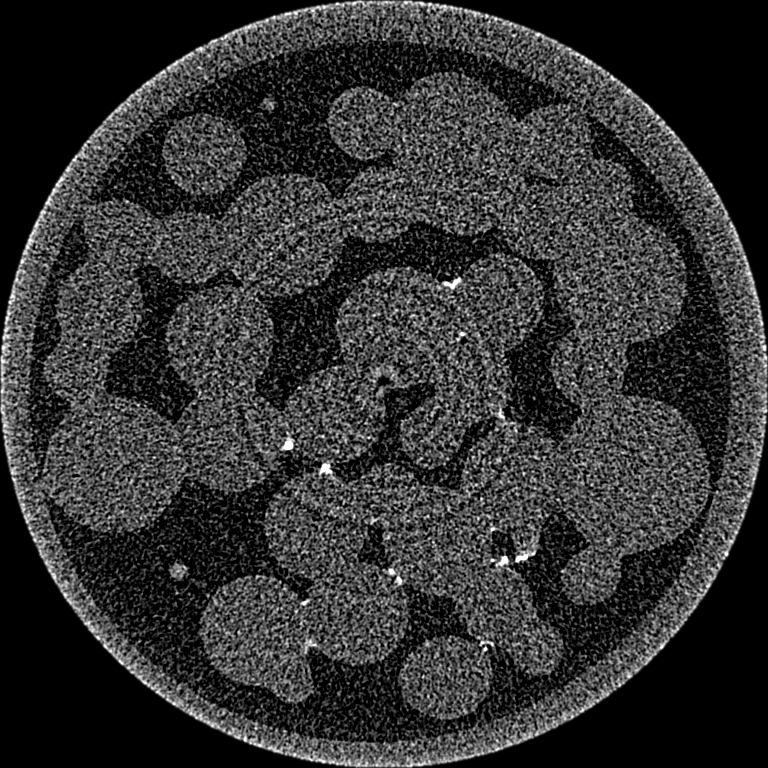}};
  \spy on (b2) in node at (0.745,-0.845);
  \node[label={[font=\footnotesize]below:AMAP-TV},inner sep=0pt] (b3) at (2,-1.1) {\includegraphics[width=0.19\textwidth]{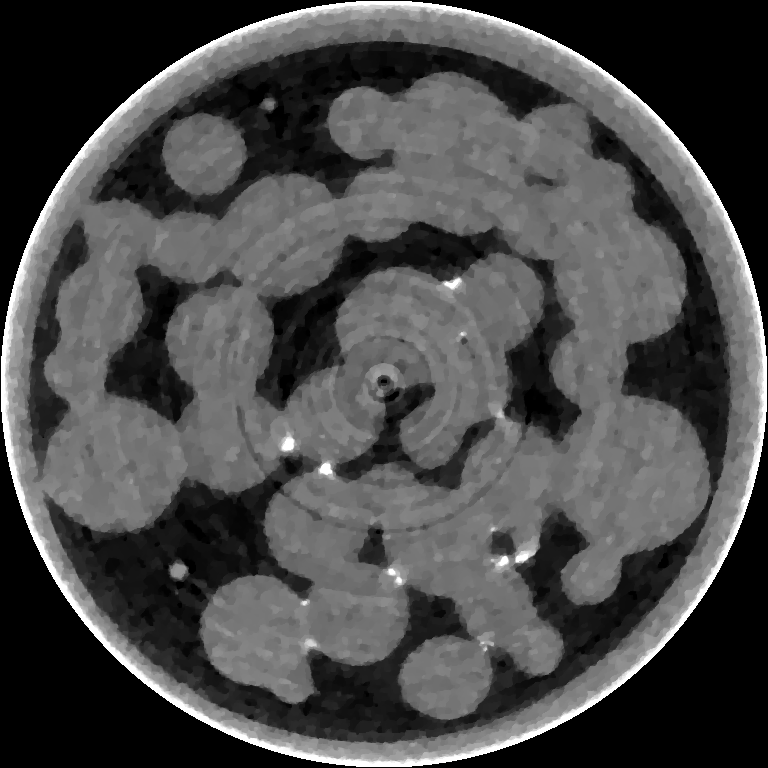}};
  \spy on (b3) in node at (1.745,-0.845);
  \node[label={[font=\footnotesize]below:JMAP-TV ($\beta=0$)},inner sep=0pt] (b4) at (3,-1.1) {\includegraphics[width=0.19\textwidth]{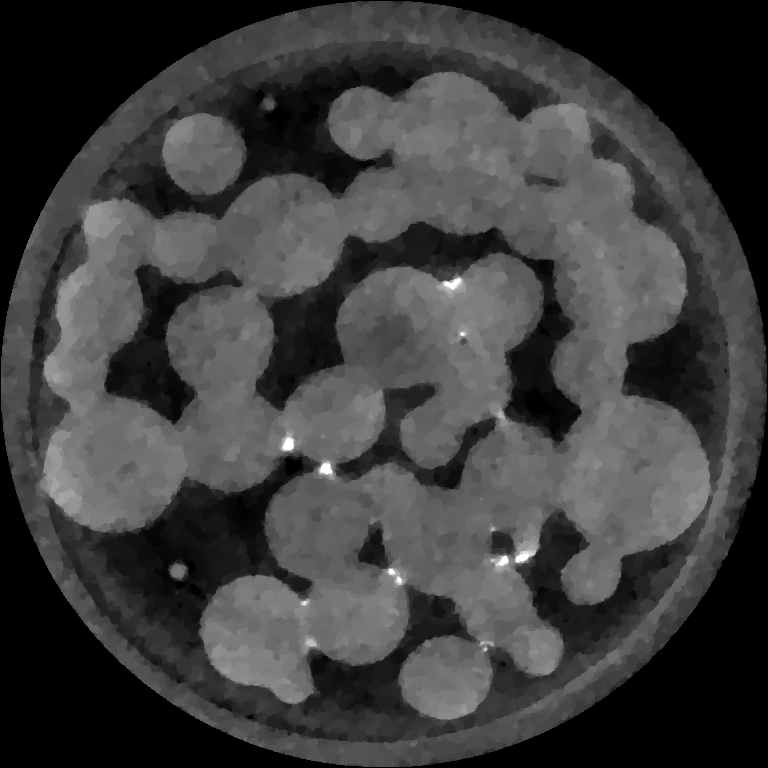}};
  \spy on (b4) in node at (2.745,-0.845);
  \node[label={[font=\footnotesize]below:JMAP-TV ($\beta=200$)},inner sep=0pt] (b5) at (4,-1.1) {\includegraphics[width=0.19\textwidth]{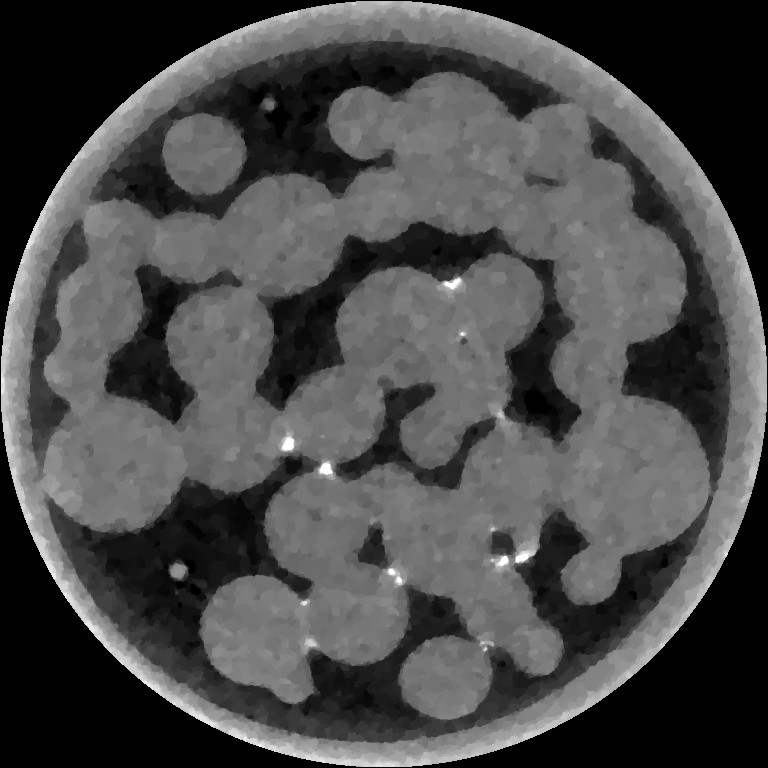}};
  \spy on (b5) in node at (3.745,-0.845);
\end{tikzpicture}

\caption{Reconstructions of real tomographic measurements. The display
  range for the images is 0 to 10 cm$^{-1}$. The reconstructions using
  the TV-prior were obtained with $\gamma = 0.01$. The number of
  iterations were 50 for reconstructions without TV prior and 1,000
  with TV prior. The insets are blow-ups of the reconstructions at the
  isocenter.}
\label{recnAPS}
\end{figure*}

We now evaluate the performance of the proposed model based on real
measurement data provided by the Advanced Photon Source (APS) facility
operated by Argonne National Laboratory (USA). The data set provides
tomographic measurements of a sample of glass beads with some dried
potassium from $p=900$ projection angles between 0$^\circ$ and
180$^\circ$ in a parallel beam geometry and with a $600\times 960$
pixel detector array. In this experiment, we will consider only a
2D reconstruction of the center slice (slice 300) so we take
$r=960$. \remove{The exposure time for each projection was only 6 ms, and the
measurements are very noisy with most photon counts in the range 0-20. We
used a square grid with side length 0.3053 cm and $768 \times 768$ pixels 
for the reconstructions. The data set contains a total of 20
flat-field measurements, but 8 of these appear to be corrupted, so we
use $s=12$ flat-field measurements.  Our reconstructions are shown in
the Fig.~\ref{recnAPS}.}\add{ The energy of the X-ray source was
33.27 keV, and the photon flux per pixel in each projection was
approximately 1200 photons/s. With an exposure time of only 6 ms, that
amounts to pixelwise photon counts in the range 0-20 per
projection. Out of a total of 20 flat-field measurements collected
before and after the experiment, 8 appear to
be corrupted, so we used $s=12$ flat-field measurements for our
reconstructions. Moreover, we used a square grid with side length 0.3053 cm and $768 \times 768$ pixels for the reconstructions. Our reconstructions are shown in the Fig.~\ref{recnAPS}.
}

Without the TV-prior on the attenuation image, the reconstructions are
quite noisy because of the low SNR. The FBP reconstruction and the
AMAP reconstruction both have ring artifacts which heavily distort the
reconstruction. The P-FBP reconstruction does not have noticeable ring
artifacts, but the reconstruction is quite noisy. Thus, to reduce
noise, we smoothed the FBP and P-FBP reconstructions using a Gaussian
filter with standard deviation $1.0$, and although this help, the
resulting images are still somewhat noisy compared to the other
reconstructions. The JMAP reconstruction with the UP prior
($\beta = 0$) has no noticeable ring artifacts, but it has a
significant amount of noise.  This is especially noticeable near the
circular boundary of the object, and it may be because of flat-field
estimation errors. Indeed, using the FE prior with $\beta=200$ yields
a reconstruction that is somewhat improved near the outer
circles. Notice that the JMAP reconstructions do not have such a
``hole'' in the middle like the FBP, P-FBP, and AMAP
reconstruction. Finally, including the TV-prior on $u$ results in the
AMAP-TV and JMAP-TV reconstructions. These results verify the
applicability of proposed model for tomographic reconstruction based
on low-intensity measurements.

\section{Conclusion} \label{conclusion}

In X-ray computed tomography, the X-ray source intensity is typically
estimated based on a number of flat-field measurements. This
estimation introduces unavoidable errors in popular reconstruction
models such as AMAP, WLS, and FBP, and these errors lead to systematic
reconstruction errors in the form of ring artifacts. By investigating
the filtered backprojection of a line in the
sinogram, we have demonstrated that such systematic errors introduce
structural changes in the reconstruction in the form of a ring. Based
on the statistics of X-ray measurements, our analysis shows an inverse
relationship between severity of ring artifacts and the source
intensity. Therefore, these systematic errors can have a significant
impact on the reconstruction quality of dose-constrained and
time-constrained problems. To mitigate this problem, we have
introduced a convex reconstruction model (JMAP) that jointly estimates
the attenuation image and the flat-field. We have also introduced a
quadratic approximation of the JMAP model, the stripe-weighted
least-squares (SWLS) model, which provides insight about the model and its
similarities with existing models.

To assess the reduction of ring artifacts in the reconstructions, we
have proposed a ``ring ratio'' error measure which quantifies the
flat-field error in the image domain. Our experimental results
indicate that the model effectively mitigates ring artifacts even for
low SNR data, not only with simulated data but also with real data
sets. In some cases, the proposed method may itself introduce
artifacts when not appropriately regularized. These artifact
essentially arise because of overfitting, and we have shown that they
can be mitigated or supressed by means of a suitable regularizing
flat-field prior. Moreover, we have shown experimentally that the JMAP
and the SWLS models have similar performance in terms of noise and
reconstruction quality.

Finally, we mention that the proposed methodology can readily be
extended to estimate a time-varying flat-field which may be useful in
applications where the flat-field does not remain stable while
acquiring the tomographic measurements and/or when the scanner acquires
projection images and flat-field images in an interleaved temporal order.

\appendices

\add{
\section{Extrema of the Radial Profile}\label{a-extrema}

The extrema of the radial profile $\tilde \mu(\rho)$, defined in \eqref{e-delta-fbp}, depend on the parameters $t_0$ and $\epsilon > 0$. To see this, we derive the critical points of $\tilde \mu(\rho)$. Setting the derivative equal to zero yields the equation
\[
\tilde \mu'(\rho) = -3\rho \left( \sigma(\sigma^2 + \rho^2)^{-5/2} +\bar \sigma(\bar \sigma^2 + \rho^2)^{-5/2} \right)= 0
\]
where $\sigma = \epsilon + i t_0$. It follows that the critical points are $\rho=0$ and any solution to the equation
\[
\sigma(\sigma^2 + \rho^2)^{-5/2} +\bar \sigma(\bar \sigma^2 + \rho^2)^{-5/2} = 0,
\]
or equivalently, $\rho = 0$ and solutions to the equation
\[
	\frac{\sigma}{\bar \sigma} = -\left(\frac{\sigma^2 + \rho^2}{\bar \sigma^2 +\rho^2}\right)^{5/2}.
\]
Taking the complex logarithm of both sides of the equation yields the equation
$
	2\angle \sigma + 2k\pi = \pi + 5 \angle (\sigma^2 +\rho^2),  k \in \mathbb{Z},
$
and hence
\begin{align}
  \label{e-angle1}
	\angle (\sigma^2 +\rho^2) = \frac{2}{5} \angle \sigma + \frac{2k-1}{5}\pi, \quad k \in \mathbb{Z}.
\end{align}
This implies that the tangent of $\angle (\sigma^2 +\rho^2)$ is equal
to  
\begin{align}
  \label{e-tan-angle1}
	\frac{2\epsilon t_0}{\rho^2+\epsilon^2-t_0^2} = \tan \left(
  \frac{2}{5} \angle \sigma + \frac{2k-1}{5}\pi\right), \quad k\in \mathbb{Z},
\end{align}
or equivalently, if we define
$c_k^{-1} = \tan \left( \frac{2}{5} \angle \sigma +
  \frac{2k-1}{5}\pi\right)$ and solve for $\rho^2$, we get
$ \rho^2 = 2\epsilon t_0 c_k +t_0^2 -\epsilon^2,
k\in\mathbb{Z}$. Thus, in addition to $\rho = 0$, the real roots of
the right-hand side of this equation are the critical points of $\tilde \mu(\rho)$, and
hence we may limit our attention to $k \in \mathbb{Z}$ for which
$2\epsilon t_0 c_k +t_0^2 -\epsilon^2 \geq 0$.

In order to find the extrema of $\tilde \mu(\rho)$, we now rewrite \eqref{e-delta-fbp}
as
\[
	\tilde \mu(\rho) = \frac{1}{4\pi} \frac{|\sigma|}{|\sigma^2+\rho^2|^{3/2}} \cos(\angle \sigma - \angle (\sigma^2 + \rho^2)).
\]
At a nonzero critical point $\rho_k \neq 0$, the angle $\angle
(\sigma^2 + \rho_k^2)$ is given by \eqref{e-angle1}, and it follows from \eqref{e-tan-angle1} that
\[
	|\sigma^2+ \rho_k^2| = 2\epsilon |t_0|\left(c_k^2 + 1\right)^{1/2}.
\]
This allows us to express the extrema associated with $\rho_k$ as
\[
	\tilde \mu( \rho_k) =
        \frac{(\epsilon^2+t_0^2)^{1/2}}{4\pi(1+c_k^2)^{3/4}(2\epsilon
          |t_0|)^{3/2} } \cos\left(\frac{2}{5} \angle\sigma +
          \frac{2k-1}{5}\pi\right), 
\]
and it immediately follows that for $|t_0| \gg  \epsilon$, the extrema are approximately inversely proportional to $\sqrt{\epsilon^3 |t_0|}$.


\section{Interpretation of Flat-field Estimate} \label{flat-field-intr}

The $i$th element of flat-field estimate $\hat v$, defined in
\eqref{e-map-vmap}, is given by
\begin{align}\label{vmap-ele}
\hat v_i(u) = \frac{\ones^T f_i + \ones^T y_i + \alpha_i - 1}{d_i(u)}
\end{align}
where $f_i \in \reals^{s}$, $y_i \in \reals^{p}$, $d_i(u) = s +
\tau_i(u) + \beta_i$, and $\tau_i(u) = \sum_{j=1}^p
\exp(-e_i^TA_ju)$. This expression can be reformulated as
\begin{align}
\hat  v_i(u) &= \frac{s}{d_i(u)} \frac{\ones^T f_i}{s} + \frac{\tau_i(u)}{d_i(u)} \frac{\ones^T y_i}{\tau_i(u)} + 
  \frac{\beta_i}{d_i(u)} \frac{\alpha_i-1}{\beta_i} \nonumber \\
  		   &=  \frac{s}{d_i(u)} (\vml)_i + \frac{\tau_i(u)}{d_i(u)} (\vmly)_i +
  \frac{\beta_i}{d_i(u)} \hat v_{\mathrm{pr}}(\alpha_i,\beta_i) \label{vmap-ele-expr}
\end{align}
where the ML estimate $\vml$ is defined in \eqref{empmeanest}, the
estimate $\vmly(\hat u)$ is defined in \eqref{e-mlu-aml}, and \[\hat
  v_{\mathrm{pr}}(\alpha,\beta) = \diag(\beta)^{-1}(\alpha - \ones)\]
is the mean of the Gamma prior. It follows from the definition
\eqref{e-c-du}, \ie, $d_i(u) = s + \tau_i(u) + \beta_i$, that 
\[ \frac{s}{d_i(u)} + \frac{\tau_i(u)}{d_i(u)} + \frac{\beta_i}{d_i(u)} = 1\]
and hence $\hat v_i(u)$ is a convex combination of three estimates.
Thus, the full flat-field vector $\hat v(u)$ can be expressed as 
\begin{align*}
  \hat v(u) = \diag(\theta_1) \vml + \diag(\theta_{2}) \vmly(\hat u) + \diag(\theta_3) \hat v_{\mathrm{pr}}(\alpha,\beta)
\end{align*}
where $\theta_1 = \diag(d(u))^{-1} s\ones$, $\theta_2 = \diag(d(u))^{-1} \tau(u)$, and $\theta_3 = \diag(d(u))^{-1} \beta$ with $\theta_1+\theta_2+\theta_3 = \ones$.
 
\section{Type-II ML Estimation of Hyperparameters} \label{type-II-estimates}
The marginal probability of $f_{i1},\ldots,f_{is}$ given the
hyperparameters $\alpha_i$ and $\beta_i$ can be computed analytically and is given by
\begin{multline}
		\prob(f_{i1},\ldots,f_{is} \mid \alpha_i, \beta_i) \\
			= \int_{0}^{\infty} \prob(
                        f_{i1},\ldots,f_{is}\mid v_i ) \prob(v_i \mid \alpha_i, \beta_i) d v_i \\
				= \frac{\Gamma(k_i + \alpha_i)}{\left(\prod_{k=1}^s f_{ik}!\right) \Gamma(\alpha_i) \, s^{k_i}} \left(\frac{\beta_i}{s+\beta_i}\right)^{\alpha_i} \left( \frac{s}{s+\beta_i}\right)^{k_i}
\end{multline}
where $k_i = \sum_{k=1}^s f_{ik}$. Here the identity
$\int_0^\infty x^b e^{-ax} \,dx = \frac{\Gamma(b+1)}{a^{b+1}}$ was
used to derive this expression. This probability distribution
resembles the negative binomial distribution, and it follows from the
first-order optimality conditions associated with \eqref{e-ml-typeII}
that $\beta_i = s\alpha_i/k_i$, or equivalently,
$\alpha_i/\beta_i = k_i/s$. This implies that the mean of the Gamma prior
is equal to the flat-field ML estimate $(\vml)_i$. Substituting the
expression for $\beta_i$ in \eqref{e-ml-typeII}, we obtain the
one-dimensional problem $\argmin_{\alpha_i} \kappa_i(\alpha_i)$ where
\begin{align*}
\kappa_i(\alpha_i) = - \log \frac{\Gamma(k_i + \alpha_i)}{\Gamma(\alpha_i)} - \alpha_i \log\frac{\alpha_i}{\alpha_i + k_i} - k_i \log\frac{k_i}{\alpha_i + k_i}.
\end{align*}
The derivative of $\kappa_i(\alpha_i)$ is
\begin{align*}
\kappa_i'(\alpha_i) &= - \left[ \digamma(k_i+\alpha_i) - \digamma(\alpha_i)  - \log \left( 1 + \frac{k_i}{\alpha_i} \right) \right] \nonumber \\
		&= - \sum_{l=0}^{k_i-1} \frac{1}{\alpha_i + l} + \log\left(1 + \frac{k_i}{\alpha_i}\right), 
\end{align*}
where $\digamma(x)$ denotes the digamma function. Similarly, the second
derivative is given by
\begin{align}
\kappa_i''(\alpha_i) = \sum_{l=0}^{k_i-1} \frac{1}{(\alpha_i + l)^2} - \frac{k_i}{\alpha_i (\alpha_i + k_i)}
\end{align}
where the summation satisfies the inequality
\begin{align}
\sum_{l=0}^{k_i-1} \frac{1}{(\alpha_i + l)^2} = \sum_{n=\alpha_i}^{\alpha_i + k_i - 1} \frac{1}{n^2} & > \int_{\alpha_i}^{\alpha_i + k_i} \frac{1}{x^2} \,dx \nonumber \\
                                               & = \frac{k}{\alpha_i (\alpha_i + k_i)}
\end{align}
for $\alpha_i > 0$. This shows that $\kappa_i''(\alpha_i)>0$ for
$\alpha > 0$, and hence $\kappa_i$ is convex on the positive
real line. Moreover, since $\kappa'(\alpha_i)$ tends to zero as as $\alpha_i$ tends to infinity, $\kappa'(\alpha_i)$ can not
have a positive zero. Consequently, the resulting flat-field Gamma prior has
zero variance (\ie, $\alpha_i/\beta_i^2$ tends to zeros for $\alpha_i
\rightarrow \infty$ since $\beta_i = s\alpha_i/k_i$) and its mean is
equal to the empirical flat-field estimate, \ie, $\alpha_i/\beta_i =(\vml)_i $.
}

\section*{Acknowledgment}
This work was supported in part by the European Research Council under Grant
No.\ 291405 (HD-Tomo) and in part by NIH R01
Grant No.\ CA158446. The contents of this article are solely the
responsibility of the authors and do not necessarily represent the
official views of the National Institutes of Health. Portions of this work were performed at
GeoSoilEnviroCARS (The University of Chicago, Sector 13), Advanced
Photon Source (APS), Argonne National Laboratory. GeoSoilEnviroCARS is
supported by the National Science Foundation -- Earth Sciences
(EAR-1128799) and Department of Energy -- GeoSciences
(DE-FG02-94ER14466). This research used resources of the Advanced
Photon Source, a U.S.\ Department of Energy (DOE) Office of Science
User Facility operated for the DOE Office of Science by Argonne
National Laboratory under Contract No. DE-AC02-06CH11357.

\ifCLASSOPTIONcaptionsoff
  \newpage
\fi
\IEEEtriggeratref{27}


%



\bibliographystyle{IEEEtran}
\bibliography{IEEEabrv,flatfield}

%



\begin{IEEEbiography}[{\includegraphics[width=1in,height=1.25in,clip,keepaspectratio]{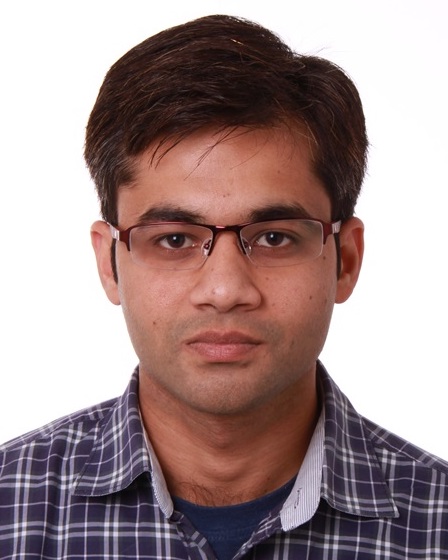}}]{Hari
    Om Aggrawal} received his M.Tech.\ in Electrical Engineering from the
  Indian Institute of Technology Kanpur, India, in 2011, and he is
  presently pursuing a PhD degree at the Technical University of
  Denmark in the Section for Scientific Computing at the Department of
  Applied Mathematics and Computer Science. His interests include
  image reconstruction models, methods for tomographic imaging, and image registration.
\end{IEEEbiography}

\begin{IEEEbiography}[{\includegraphics[width=1in,height=1.25in,clip,keepaspectratio]{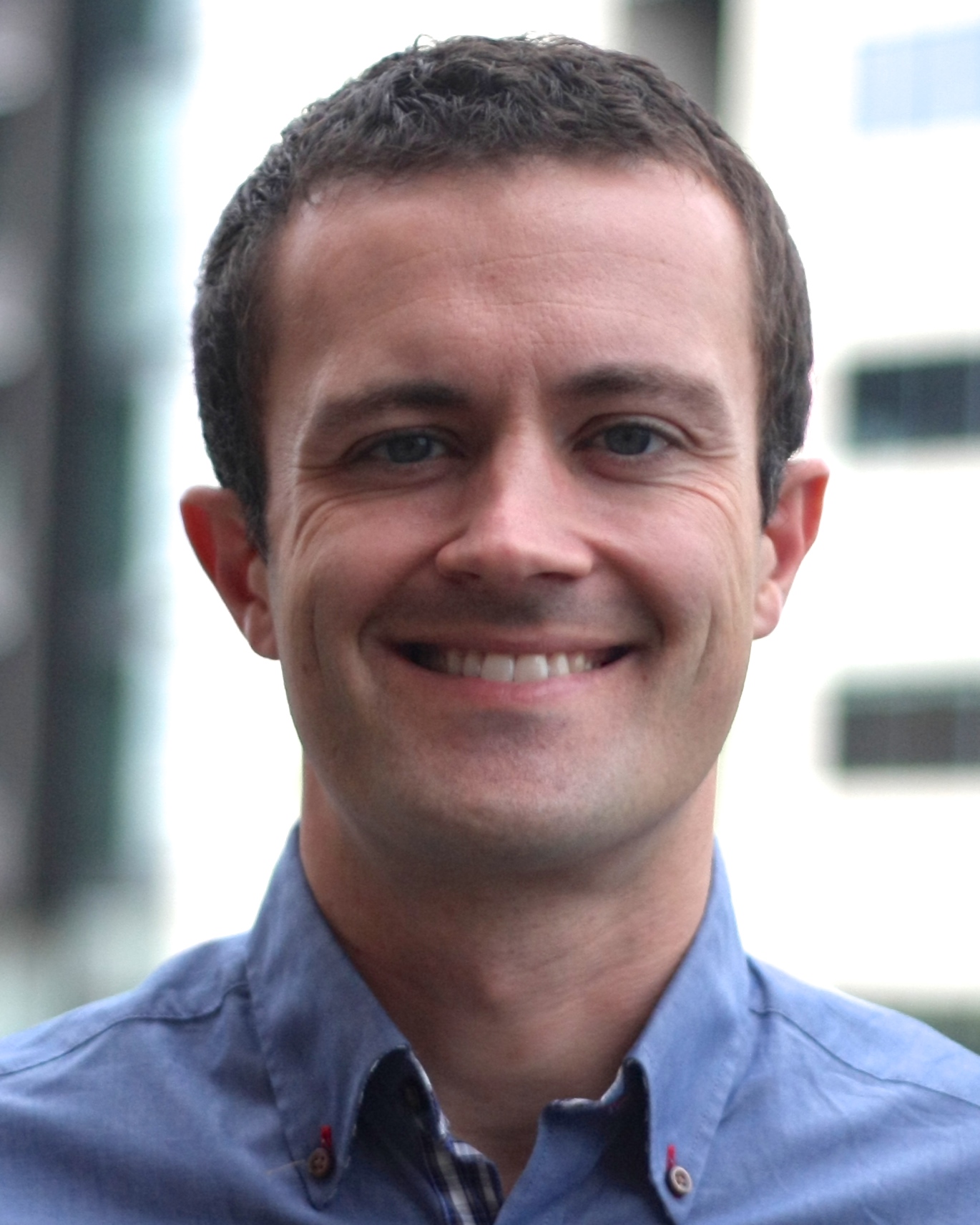}}]{Martin
    S. Andersen} received his M.S.\ in Electrical Engineering from
  Aalborg University, Denmark, in 2006 and his PhD in Electrical
  Engineering from the University of California, Los Angeles, in
  2011. After receiveing his degree, he was a postdoc in the Division
  of Automatic Control at Link{\"o}ping University, Sweden, and at the
  Technical University of Denmark. He is currently an Associate
  Professor at the Technical University of Denmark in the Section for
  Scientific Computing at the Department of Applied Mathematics and
  Computer Science. His research interests include optimization,
  numerical methods, signal and image processing, and systems and
  control.
\end{IEEEbiography}

\begin{IEEEbiography}[{\includegraphics[width=1in,height=1.25in,clip,keepaspectratio]{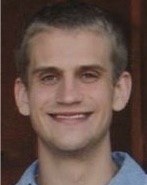}}]{Sean
    D. Rose} was born in Cincinnati, OH. He received a B.S.\ degree 
in physics and biochemistry from The Ohio State University in 2013.
He began his graduate studies in Medical Physics at the University of Chicago 
in the summer of 2013 in the lab of Dr.\ Xiaochuan Pan where he continues
to work on the development and optimization of iterative methods for tomographic
image reconstruction. His research interests are in optimization based 
image reconstruction for tomographic imaging modalities, objective assessment 
of image quality and its application to parameter selection, and large-scale optimization. 
\end{IEEEbiography}

\begin{IEEEbiography}[{\includegraphics[width=1in,height=1.25in,clip,keepaspectratio]{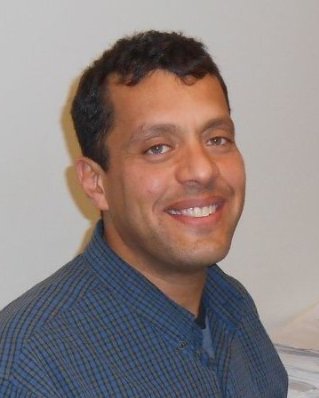}}]{Emil
    Y. Sidky} has held academic positions in Physics at the University of 
Copenhagen and Kansas State University. He then switched to Medical
Physics joining the University of Chicago in 2001, where he is now a Research Associate 
Professor. His current interests are image reconstruction in X-ray computed tomography
and digital breast tomosynthesis, large-scale optimization, and 
objective assessment of image quality.
\end{IEEEbiography}


\vfill


\end{document}